\input amstex
\documentstyle{amsppt}
\pageno=1
\magnification1200
\catcode`\@=11
\def\logo@{}
\catcode`\@=\active
\NoBlackBoxes
\vsize=23.5truecm
\hsize=16.5truecm

\def\d{d\!@!@!@!@!@!{}^{@!@!\text{\rm--}}\!}

\def\leg{\;\dot{\le}\;}

\def\crp{\overline{\Bbb R}_+}
\def\crm{\overline{\Bbb R}_-}

\def\rnp{{\Bbb R}^n_+}

\def\rnpm{\Bbb R^n_\pm}
\def\crnp{\overline{\Bbb R}^n_+}

\def\Rn{\Bbb R^n}
\def\ang#1{\langle {#1} \rangle}

\def\ttilde{\overset{\,\approx}\to}

\def\Zfrac{\tsize\frac1{\raise 1pt\hbox{$\scriptstyle z$}}}

\def\rp{ \Bbb R_+}
\def\rmi{ \Bbb R_-}
\define\tr{\operatorname{tr}}
\define\op{\operatorname{OP}}

\define\Tr{\operatorname{Tr}}

\define\stimes{\!\times\!}

\document

\topmatter
\title On the logarithm component in trace defect formulas 
\endtitle
\author
Gerd Grubb
\endauthor
\affil
{Copenhagen Univ\. Math\. Dept\.,
Universitetsparken 5, DK-2100 Copenhagen, Denmark.
E-mail {\tt grubb\@math.ku.dk}}\endaffil
\abstract 
In asymptotic expansions of resolvent traces $\Tr(A(P-\lambda
)^{-1})$ for classical pseudodifferential operators on closed
manifolds, the coefficient $C_0(A,P)$ of $(-\lambda )^{-1}$ is of
special interest, since it is the first coefficient containing
nonlocal elements from $A$; on the other hand if $A=I$ and $P=D^*D$
it gives part of the index of $D$. $C_0(A,P)$ also equals the
zeta function value at $0$ when $P$ is invertible.
$C_0(A,P)$ is a trace modulo local
terms, since $C_0(A,P)-C_0(A,P')$ and $C_0([A,A'],P)$ are local. By
use of complex powers $P^s$ (or similar holomorphic families
of order $s$), Okikiolu, Kontsevich and Vishik, Melrose and Nistor showed
formulas for these trace defects in terms of residues of operators
defined from $A$, $A'$, $\log P$ and $\log P'$.

The present paper has two purposes: One is to show how the trace
defect formulas can be obtained from the resolvents in a simple way
without use of the 
complex powers of $P$ as in the original proofs. We here also give a simple
direct proof of a recent residue formula of Scott for $C_0(I,P)$. The
other purpose is to establish trace defect residue formulas for
operators on manifolds with 
boundary, where complex powers are not easily accessible; we do this
using only resolvents. We also
generalize Scott's formula to boundary problems.
\endabstract
\rightheadtext{The logarithm component}
\endtopmatter

\subhead Introduction \endsubhead

Consider a classical pseudodifferential operator ($\psi $do) $A$ of
order $\sigma $ on an $n$-dimensional smooth
compact boundaryless manifold $X$. When $P$ denotes an auxiliary
elliptic $\psi $do 
of order $m>0$ and, say, positive, one can study the generalized zeta
funcion $\zeta (A,P,s)$ defined as the meromorphic extension of
$\Tr(AP^{-s})$ to the complex plane, where the complex powers $P^{-s}$
are defined from the resolvent $(P-\lambda )^{-1}$ as in Seeley
\cite{S}. It is well-known that $\zeta (A,P,s)$ 
has a Laurent expansion at $s=0$,
$$
\zeta (A,P,s)\sim C_{-1}(A,P)s^{-1}+C_0(A,P)+\sum_{l\ge
1}C_l(A,P)s^l,\tag 0.1
$$
where $m\, C_{-1}(A,P)$ equals the noncommutative residue 
$\operatorname{res}A$ (Wodzicki \cite{W}, Guil\-le\-min \cite{Gu}), and
$C_0(A,P)$ equals the canonical trace $\operatorname{TR} A$ in
particular cases 
(Kontsevich and Vishik \cite{KV}, Lesch \cite{L}, recent extension in 
Grubb \cite{G2}).

The coefficient $C_0(A,P)$ is not in general independent of $P$, but
then it is viewed as a ``regularized trace'' (Melrose and Nistor
\cite{MN}) or a ``weighted trace'' (Cardona, Ducourtioux, Magnot and Paycha
\cite{CDMP}, \cite{CDP}). In general it satisfies the {\it trace
defect formulas}
$$
\gather
C_0(A,P)-C_0(A,P')=-\tfrac 1m \operatorname{res}(A(\log P-\log P')),\tag0.2
\\
C_0([A,A'],P)=-\tfrac 1m \operatorname{res}(A[A',\log P]),\tag0.3
\endgather
$$ 
shown by Okikiolu \cite{O} and \cite{KV}, resp.\ \cite{MN}, by use
essentially of the holomorphic family $P^{-s}$ and the fact that its
derivative at 0 is $-\log P$.

For a compact manifold $X$ with boundary $\partial X=X'$, the
situation is somewhat different. A calculus that contains
differential elliptic boundary value problems and their solution
operators and has full pseudodifferential composition rules is the
calculus of Boutet de Monvel \cite{B}; we consider an operator
$A=P_++G$ lying there. Here $P$ is a $\psi $do defined on a larger
boundaryless manifold $\widetilde X$ in which $X$ is imbedded, such
that $P$ satisfies the transmission condition at $X'$ (in particular,
it is of integer order), and $G$ is a singular Green operator
(smoothing in the interior, but important near the boundary).

Even the simplest auxiliary operator $P _{1,\operatorname{D}}$
with $P_1$ equal to the Laplace operator, the D indicating Dirichlet
condition, 
does {\it not} have its complex powers in the Boutet de
Monvel calculus, so the ingredients in the zeta function $\Tr
(AP_{1,\operatorname{D}}^{-s})$ are not easily accessible.
Nevertheless, by passing via the resolvent family
$A(P_{1,\operatorname{D}}-\lambda )^{-1}$, we managed to show in a
joint work with Schrohe \cite{GSc2}, that
$C_0(A,P_{1,\operatorname{D}})-C_0(A,P_{2,D})$ and
$C_0([A,A'],P_{1,\operatorname{D}})$ are {\it local}, and to pinpoint
the nonlocal content of $C_0(A,P_{1,\operatorname{D}})$ modulo local terms.
The question of possible generalizations of the formulas (0.2)--(0.3)
remained open then.

\medskip

In the present paper we show for the boundaryless case how the formulas
(0.2)--(0.3) can be derived
directly from the knowledge of the resolvent (Section 2). The
crucial fact is that the constant comes from a strictly homogeneous
term in the symbol of \linebreak$A((P-\lambda )^{-1}-(P'-\lambda )^{-1})$
resp.\ $A[A',(P-\lambda )^{-1}]$ which is integrable at $\xi =0$ (and
is $O(\lambda ^{-2})$ for $|\lambda |\to\infty $ when $\xi \ne 0$). The
operator $\log P$ appears simply 
because $\log \lambda $ has a jump of $2\pi i$ at the negative real
axis; there is no need to construct the $P^{-s}$.

Before this we give (in Section 1) a similarly simple proof of the formula
shown recently by Scott \cite{Sco}:
$$
C_{0}(I,P)=-\tfrac1m\operatorname{res}(\log P),\tag0.4
$$
from which he draws consequences on multiplicative properties; here
$C_0(I,P)$ equals $\zeta (P,0)$ $+$ the nullity of $P$.
Scott's proof of (0.4) is based on calculations inspired from
\cite{O}, going via results for $P^{-s}$. In fact, finding the
direct proof of (0.4) in terms of the resolvent was the starting
point for our present paper.

Next, we discuss possible generalizations of the formulas
(0.2)--(0.3) to the situation with boundary. Here we replace the
family $(P_{1,\operatorname{D}}-\lambda )^{-1}$ used in \cite{GSc1,2} by
its $ \psi $do part $(P_1-\lambda )^{-1}_+$ (which corresponds to
replacing   $P_{1,\operatorname{D}}^{-s}$ by  $(P_1^{-s})_+$, another
family which equals the identity for $s=0$); this spares us for the
technicalities involved in working with a boundary
condition. On the other hand we
allow general higher order choices of $P_1$, where \cite{GSc1,2}
considered the second-order case (which provides simple roots
in the detailed construction of the resolvent symbol). To handle
general choices of $P_1$, we base the study on the relatively crude
methods from the book \cite{G1}.

In Section 3 we show that (0.2) does generalize in a natural way, since
$(\log P_1-\log P_2)_+$ is a zero order $\psi $do having the
transmission property:
$$
C_0(A,P_{1,+})-C_0(A,P_{2,+})=-\tfrac 1m \operatorname{res}(A(\log
P_1-\log P_2)_+).\tag 0.5
$$ Here the residue definition of Fedosov, Golse, Leichtnam and
Schrohe \cite{FGLS} is used. 

In Section 4 we consider generalizations of (0.3), for two operators
$A=P_++G$, \linebreak $A'=P'_++G'$ of orders $\sigma $ and $\sigma '$, and
normal order 0. The leftover
terms (singular Green type terms) in commutators $[A',(\log P_1)_+]$ are 
not in the calculus and have not (yet) been covered 
by residue formulas, so we cannot extend (0.3) directly. 
However, considering $A[A',(P_1-\lambda )_+^{-1}]$, we show that the
normal trace $\Cal S_\lambda $ of its singular Green operator part
$\Cal G_\lambda $ 
is a $\psi $do on $X'$ with sufficiently good symbol estimates to
allow integration against $\log \lambda $, leading to
a classical $\psi $do $S$ on $X'$ such that
$$
C_0([A,A']P_{1,+})=-\tfrac1m\operatorname{res}_X((P[P',\log P_1
])_+)-\tfrac1m\operatorname{res}_{X'}(S).\tag 0.6
$$
\medskip

Finally in Section 5, we show a certain generalization of (0.4) to
elliptic pseudodifferential boundary problems $(P_++G)_T$ as
considered in \cite{G1}.

\subhead 1. On the residue of logarithm formula \endsubhead

Let $P$ be an elliptic pseudodifferential operator of order $m \in
\Bbb R_+$ on a closed (i.e., compact boundaryless)
manifold of dimension $n$, such that the
principal symbol has no eigenvalues on $\Bbb R_-$. We can assume
that $P$ has no eigenvalues on $\Bbb R_-$ (by a small rotation
if needed). Then we can
define the resolvent $Q_\lambda =(P-\lambda )^{-1}$ in a
sector $V$ around $\Bbb R_-$. The complex powers and the logarithm are
defined by functional calculus:
$$\aligned
P^{-s}&=\tfrac i{2\pi }\int_{\Cal C}\lambda
^{-s}(P-\lambda )^{-1}\,d\lambda \text{ for }
\operatorname{Re}s>0,\quad P^{k-s}=P^kP^{-s};\\ 
\log P&=\lim_{s\to 0}\tfrac i{2\pi }\int_{\Cal C}\lambda
^{-s}\log\lambda \,(P-\lambda )^{-1}\,d\lambda ,
\endaligned\tag1.1
$$
with integrations on a curve $\Cal C$ in $\Bbb C\setminus \overline{\Bbb
R}_-$ going around the nonzero spectrum of $P$ in the positive
direction; hereby $P^{-s}$ and $\log P$ are taken to be 0 on $\ker
P$. 

It is well-known 
that $\Tr P^{-s}$ extends meromorphically to $\Bbb C$ as the zeta function
$\zeta(P,s)$ (Seeley \cite{S}); it
is regular at $s=0$.
 It is
also known 
that the noncommutative residue can be defined for $\log P$ (Okikiolu
\cite{O}, Lesch \cite{L}). The value at $s=0$ was recently identified
by Scott (\cite{Sco}) with a residue: 
$$
\zeta(P,0)=-\tfrac1m \operatorname{res}(\log P)\tag1.2
$$
(this is the formula if $\ker P=0$; also nonzero cases are considered).
His method is based on an analysis of the symbol of $P^{-s}$
inspired from \cite{O}. We shall show below how the formula can be
proved directly from the knowledge of the resolvent.

We assume $m >n$
for convenience. (Otherwise, one can consider
$(P-\lambda )^{-N}$ for large $N$, where the local formulas however
boil down to the same calculation, as indicated in a general
situation in Remark 3.12 below.)
Then $Q_\lambda $ is trace-class, and its kernel (calculated in local
coordinates) has an asymptotic expansion for
$\lambda \to\infty $ in $V$, leading to a trace expansion by
integration of the fiber trace in $x$:
$$
\aligned
K( Q_\lambda ,x,x)
&\sim
\sum_{ j\ge 0  }   c_{ j}(x)(-\lambda ) ^{\frac{n -j}m  -1}+ 
\sum_{ k\ge 1}\bigl(   c'_{ k}(x)\log (-\lambda ) +  c''_{
k}(x)\bigr)(-\lambda ) ^{ -k-1},\\
\Tr Q_\lambda
&\sim
\sum_{ j\ge 0  }   c_{ j}(-\lambda ) ^{\frac{n -j}m  -1}+ 
\sum_{ k\ge 1}\bigl(   c'_{ k}\log (-\lambda ) +  c''_{
k}\bigr)(-\lambda ) ^{ -k-1}.
\endaligned\tag1.3
$$
This was first shown by Agranovich \cite{A} (with reference to the heat trace
formulation of Duistermaat-Guillemin \cite{DG} and the complex power
formulation of Seeley \cite{S}); proof details can also be found in
Grubb and Seeley \cite{GS1} for the
case where $m $ is integer and in Loya \cite{Lo}, Grubb and Hansen
\cite{GH} for the 
general case.  In fact, the meromorphic structure of $\zeta (P,s)$
and the asymptotic expansion of $\Tr Q_\lambda $ can be deduced from
one another (as accounted for e.g.\ in \cite{GS2}).
In particular,  we can define
$$
C_0(P)=c_n=\int_{X}\tr c_n(x)\,dx;\;\text{ then }C_0(P)=\zeta (P,0)+\nu
_0,\tag1.4
$$ 
where $\nu _0$ is the algebraic multiplicity of 0 as an eigenvalue of
$P$. For, $\nu _0$ equals the rank of the eigenprojection $\Pi
_0=\tfrac i{2\pi }\int_{|\lambda |=\varepsilon } 
(P-\lambda )^{-1}\,d\lambda $, cf.\ Kato \cite{K,
Sect.\ III 6.8}.

We shall base our study of $C_0(P)$ on the resolvent information, and
will now recall an
elementary deduction of the kernel expansion down to $O(|\lambda |
^{-2+\varepsilon })$.
In local coordinates, the symbol $q(x,\xi ,\lambda )$ of $Q_\lambda $ has an expansion in
quasi-homogeneous terms $
q(x,\xi ,\lambda )\sim\sum_{j\ge 0}q_{-m-j}(x,\xi ,\lambda )$,
where $q_{-m}=(p_m-\lambda )^{-1}$, and $q_{-m-j}$ for each $j\ge 1$
is a finite sum of terms with the structure 
$$
f(x,\xi ,\lambda )=g_1\,q_{-m}^{\nu _1}\,g_2\,q_{-m}^{\nu _2}\dots g_M\,q_{-m}^{\nu _M}\,g_{M+1};\tag1.5
$$
here the $\nu _k$ are integers $\ge 1$ and the $g_k(x,\xi )$ are
$\psi $do symbols independent of $\lambda $ and homogeneous of degree
$r_k$ for $|\xi |\ge 1$. The index sums $r=\sum_{1\le k\le
M+1}r_k$ and $\nu 
=\sum_{1\le k\le M}\nu _k$ satisfy$$
2\le \nu \le 2j+1,\quad r=-j+(\nu -1)m.\tag1.6
$$
This is seen by working out the symbol construction in \cite{S} in detail
(more information and references in \cite{G1, Rem.\ 3.3.7}).
We indicate strictly homogeneous versions (the extensions by
homogeneity into the region $|\xi |\le 1$) by an upper index
$h$; the $q^h_{-m-j}$ satisfy
$$
q^h_{-m-j}(x,t\xi ,t^m\lambda )=t^{-m-j}q^h_{-m-j}(x,\xi
,\lambda )\text{ for }t>0, \text{ all }\xi \ne 0.\tag1.7 
$$

Note that in (1.5), $f^h$ is $O(|\xi |^{r})$ at $\xi =0$ (for
$\lambda \ne 0$), hence
integrable in $\xi $ at $\xi =0$
if $r>-n$. Then in view of (1.6), $q^h_{-m-j}$ is integrable at $\xi
=0$ when $j<n+m$ and $\lambda \ne 0$ (this is clear for $j=0$, and for $j\ge 1$, the least integrable contributions
are those with 
$\nu =2$). In particular, $q^h_{-m-n}$ is continuous in $\xi $.

The diagonal kernel $K(Q_\lambda ,x,x)$ defined from $q$ equals $\int_{\Bbb
R^n}q(x,\xi ,\lambda )\,\d\xi $ (where  $\d$ stands for $ (2\pi )^{-n}d$).

\proclaim{Lemma 1.1} $q$ has an expansion in
strictly homogeneous terms plus a remainder:$$
q(x,\xi ,\lambda )=\sum_{0\le j<m+n}q^h_{-m-j}(x,\xi ,\lambda
)+q'_{-2m-n}(x,\xi ,\lambda ),\tag1.8
$$
where the $q^h_{-m-j}$ ($j<m+n$) and $q'_{-2m-n}$ are integrable in
$\xi $, and 
$\int q'_{-2m-n}\,\d\xi =O(|\lambda |^{-2+\varepsilon })$, any
$\varepsilon >0$. 
Consequently, $K(Q_\lambda ,x,x) $ has the expansion$$
\aligned
K( Q_\lambda ,x,x)
&=
\sum_{ j\ge 0  }   c_{ j}(x)(-\lambda ) ^{\frac{n -j}m  -1}+
O(|\lambda |^{-2+\varepsilon }),\text{ where }\\
  c_j(x)&=\int_{\Bbb R^n}q^h_{-m-j}(x,\xi ,-1)\,\d\xi ,\text{ for } j<m+n.
\endaligned
\tag 1.9 
$$
\endproclaim 

\demo{Proof}
For $j=0$,$$
q_{-m}-q^h_{-m}=(p_m-\lambda )^{-1}-(p^h_m-\lambda
)^{-1}=
(p_m-\lambda )^{-1}(p^h_m-p_m)(p^h_m-\lambda )^{-1},\tag 1.10
$$
so it is supported in $|\xi |\le 1$ and $O(|\lambda |^{-2})$ there. This also
holds for $q_{-m-j}-q^h_{-m-j}$ for general $j\ge 1$ since $\nu \ge
2$ in (1.5). For $j<m+n$ the $q_{-m-j}-q^h_{-m-j}$ are integrable in $\xi $, the integrals
being $O(\lambda ^{-2})$. For the remainder $q-\sum_{j<m+n}q_{-m-n}$, write $m=m'+\delta $, $m'$ integer and $\delta
\in\,]0,1]$, and note that $j<m+n$ means
$j\le m'+n$. The symbol $q-\sum_{0\le j<m+n}q_{-m-j}$
is of order $-m-m'-n-1=-2m -n+\delta -1$ and satisfies
$$
\aligned
|q-\sum_{j<m+n}q_{-m-j}|&\le c (1+|\xi |^m+|\lambda |)^{-2}(1+|\xi
|)^{-n+\delta -1}\\
&\le c'(1+|\lambda |)^{-2+\varepsilon }(1+|\xi
|)^{-n+\delta -1-m\varepsilon },
\endaligned\tag 1.11
$$
any $\varepsilon \ge 0$. If $\delta <1$ (the case where $m$ is
noninteger), we can take $\varepsilon =0$, otherwise we take it small
positive; then the integral in $\xi $ is
$O(|\lambda |^{-2+\varepsilon })$.
This shows the statements on (1.8).

Now (1.9) follows directly by integration in $\xi $, using the
calculations
$$
\int_{\Bbb R^n}q^h_{-m-j}(x,\xi ,\lambda )\,\d\xi
=|\lambda |^{\frac{n-j}m-1}\int_{\Bbb R^n}q^h_{-m-j}(x,\eta ,\lambda
/|\lambda |)\,\d\eta .\tag1.12 
$$ 
For $\lambda \in\rmi$, they show that $
  c_j(x)=\int_{\Bbb R^n}q^h_{-m-j}(x,\xi ,-1)\,\d\xi$;
this remains valid on general rays in $V$ since $q^h_{-j-m}$ is holomorphic
in $\lambda $ (cf.\ e.g.\ \cite{GS1,
Lemma 2.3}).\qed
\enddemo 

In the case $j=n$ we get in particular, when the contributions
$c_n(x)$ are carried back to the manifold and collected:
$$
C_0(P)=  c_n=\int\tr c_n(x)\,dx,\text{ where }
c_n(x)=\int_{\Bbb R^n} q^h_{-m-n}(x,\xi ,-1)\,\d\xi .
\tag1.13
$$

Now consider the operator $\log P$, (1.1).
It is well-known that it has a symbol in local coordinates (cf.\ e.g.\ [O])$$
\operatorname{symb}(\log P)=m\log[\xi ]I+b(x,\xi ),\tag1.14
$$
where $b$ is classical of order 0, and $[\xi ]$
stands for a smooth positive function equal to $|\xi |$ for $|\xi
|\ge 1$. This symbol is
found termwise from the symbol of
$Q_\lambda =(P-\lambda )^{-1}$ by Cauchy integral formulas as in
(1.1); in particular,  
$$
b_{-n}(x,\xi )=\tfrac i{2\pi }\int_{\Cal C'}\log\lambda
\,q_{-m-n}(x,\xi ,\lambda )\, d\lambda ,\tag1.15
$$
where $\Cal C'$ is a closed curve in $\Bbb C\setminus\crm$ encircling
the eigenvalues of 
$p_m$.
According the the definition of noncommutative residues of operators
with log-polyhomogeneous symbols (\cite{O}, \cite{L}), $$
\operatorname{res}(\log P)=\int_{X}\int_{|\xi |=1}\tr b_{-n}(x,\xi
)\,\d S(\xi )dx\tag1.16
$$
(where the integral is known to have an invariant meaning).
We want to show that this number equals $-mC_0(P)$.
This will be based on a simple lemma.

\proclaim{Lemma 1.2}Let $f(\lambda )$ be meromorphic on $\Bbb C$ and
$O(\lambda ^{-1-\varepsilon })$ for $|\lambda |\to\infty $ (some
$\varepsilon >0$), with
poles lying in a bounded subset of $\Bbb C\setminus\crm$. Let
$\Cal C$ be a closed curve in $\Bbb C\setminus \crm$ encircling the
poles in the positive direction.
 Then$$
\tfrac 1{2\pi i}\int_{\Cal C}\log \lambda \,
f(\lambda )\,d\lambda =\int_{-\infty }^0 f(t)\,dt .\tag1.17
$$

The identity also holds if $f(\lambda )$ is holomorphic in a keyhole region
around $\crm$:$$
V_{r,\theta }=\{\lambda \in\Bbb C\mid |\lambda |< r\text{ or
}|\arg\lambda -\pi |<\theta \}\tag1.18
$$ ($r$ and $\theta >0$), and $f(\lambda )$ is
$O(\lambda ^{-1-\varepsilon })$ for $\lambda \to\infty $ in
$V_{r,\theta }$; then
$\Cal C$ should be a curve in $\Bbb C\setminus \crm$ going around
$\complement V_{r,\theta }$ in the positive direction, e.g.\ defined as the
boundary of $V_{r',\theta '}$ for some $r'\in \,]0,r[\,$, $\theta
'\in \,]0,\theta [\,$.
\endproclaim 
 \demo{Proof}
We can replace $\Cal C$ by the curve $\Cal C_1+\Cal C_2+\Cal C_3+\Cal
C_4$ in the complex plane cut-up along $\crm$, where (for a
sufficiently large $R$)
$$\alignedat3
\Cal C_1&=\{\,Re^{i\omega }\mid -\pi \le \omega \le \pi\,\}, &&\quad\
\Cal C_2&&=\{\,se^{i\pi }\mid R\ge s \ge \tfrac 1R\,\},\\
\Cal C_3&=\{\,\tfrac 1R e^{i\omega }\mid \pi \ge \omega \ge -\pi\,\}, &&\quad
\Cal C_4&&=\{\,se^{-i\pi }\mid \tfrac 1R\le s \le R\,\};
\endalignedat\tag1.19
$$
we shall let $R\to \infty $. Here $$
\aligned
\big |\int_{\Cal C_1} \log \lambda f(\lambda )\,d\lambda \big
|&=O(RR^{-1-\varepsilon }\log R)\to 0 \text{ for }R\to \infty,\\ 
\big |\int_{\Cal C_1} \log \lambda f(\lambda )\,d\lambda \big
|&=O(R^{-1 }\log R)\to 0 \text{ for }R\to \infty 
;\endaligned\tag1.20
$$
moreover, $$
\log\lambda =\log (se^{i\pi })=\log s +i\pi \text{ on } \Cal C_2
\quad
\log\lambda =\log (se^{-i\pi })=\log s -i\pi \text{ on } \Cal C_4,
$$
(the difference of the values of log $\lambda $ from above and from
below on $\Bbb R_-$ is $2\pi i$).
Then $$
\tfrac 1{2\pi i}\int_{\Cal C}\log \lambda
f(\lambda )\,d\lambda =\tfrac 1{2\pi i}\int_{-R}^{-\frac1R} 2\pi if(t)\,dt
+O(R^{-\varepsilon }\log R)=\int_{-\infty }^0f(t)\,dt.
$$

For the second statement we can instead approximate $\Cal C$ by $\Cal
C'_R=\Cal C'_1+\Cal C_2+\Cal C_3+\Cal
C_4+\Cal C'_5$, where $\Cal C_2$,  $\Cal C_3$ and $\Cal C_4$ are as
above, and
$$
\Cal C'_1=\{\,Re^{i\omega }\mid \pi-\theta ' \le \omega \le \pi\,\},
\quad
\Cal C'_5=\{\,Re^{i\omega }\mid -\pi \ge \omega \ge -\pi-\theta '\,\};
$$  
then we use that the integrals over $\Cal C'_1$ and $\Cal C'_5$ go to
0 for $R\to\infty $. 
\qed
\enddemo

At each $x$, we have the formula for $c_n(x)$ in (1.13), and the
formula with $b_{-n}$:$$
\int_{|\xi |=1}b_{-n}(x,\xi )\,\d S(\xi )=
\int_{|\xi |=1}\tfrac i{2\pi }\int_{\Cal C'}\log \lambda
\,q^h_{-m-n}(x,\xi ,\lambda )\,d\lambda 
\d S(\xi ),
\tag1.21$$
so the identification of $C_0(P)$ and
$-\frac1m\operatorname{res}(\log P)$ will be obtained if we show that for each $x$,
$$
\int_{\Bbb R^n}q^h_{-m-n}(x,\xi ,-1)\,\d\xi
 =-\tfrac1m \int_{|\xi |=1}\tfrac i{2\pi }\int_{\Cal C'}\log \lambda
\,q^h_{-m-n}(x,\xi ,\lambda )\,d\lambda 
\d S(\xi ).\tag1.22
$$

We transform the left-hand side by use of the quasi-homogeneity  (1.7).
For later reference, the calculation will be formulated in a lemma:

\proclaim{Lemma 1.3} Let $m>0$. Let $f(\xi ,t )$ be continuous
for $(\xi ,t)\in (\Bbb R^n\setminus\{0\})\times\crm$ and
quasihomogeneous there in the 
sense that $f(s\xi ,s^mt)=s^{-m-n}f(\xi ,t)$ for all $s>0$, and
integrable at $\xi = 0$ for
each $t\ne 0$. Then 
$$
\int_{\Bbb R^n}f(\xi ,-1)\,\d\xi =\tfrac1m\int_{|\xi 
|=1}\int_{-\infty }^0 f( \xi ,t)\,dt \d S(\xi  ).
\tag1.23
$$ 
\endproclaim  

\demo{Proof}
Since $|f(\xi ,-1)|=|\xi |^{-m-n}|f(\xi /|\xi |,-|\xi |^{-m})|$ is
$O(|\xi |^{-m-n})$ for $|\xi |\to \infty $, the function in the
left-hand side is integrable.
\comment
Note that for $\xi \ne 0$, $f(r\xi ,-1)r^{n-1}=r^{-m-1}f(\xi
,-r^{-m})$ where $f(\xi ,r^{-m})\to f(\xi ,0)$ for $r\to\infty $,

\endcomment
For $|\xi |=1$ we make a calculation using
the coordinate change $t=-r^{-m}$, $dt=mr^{-m-1}dr$:
$$
\int_0^\infty f(r \xi ,-1)r^{n-1}\,dr =
\int_0^\infty f(\xi ,-r^{-m})r^{-m-1}\,dr  
=\tfrac1m \int_{-\infty }^0 f( \xi ,t)\,dt,\tag1.24
$$
which gives:
$$
\int_{\Bbb R^n}f(\xi ,-1)\,\d\xi =\int_{|\eta |=1}\int_0^\infty
f(r\eta  ,-1)r^{n-1}\,dr \d S(\eta )=
\tfrac1m\int_{|\eta |=1}\int_{-\infty }^0 f( \eta ,t)\,dt \d S(\eta ),
$$
showing (1.23).
(We are using the Fubini theorem; in fact (1.24) is valid almost
everywhere with respect to $\xi \in S^{n-1}$.) 
\qed
\enddemo 

Now (1.22) follows by application of
(1.23) to $q^h_{-m-n}(x,\xi ,t )$ at each $x$ and
application of Lemma 1.2 to $\int_{-\infty }^0 q^h_{-m-n}(x,\xi ,t)\,
dt$ (the minus comes from replacing  $\frac1{2\pi i}$ by $\frac i
{2\pi }$).
Integration in $x$ of the fiber trace then gives the desired identity (1.2). 

We have shown:

\proclaim{Theorem 1.4}
$C_0(P)$ equals $-\frac1m \operatorname{res}(\log P)$, and this holds
pointwise, in 
that
$$
C_0(P)=\int_X \tr c_{n}(x)\,dx
=-\tfrac1m\operatorname{res}(\log P),
\tag1.25$$ 
where, for each $x$, in local coordinates,
$$
c_{n }(x)=\int_{\Bbb R^n}q^h_{-m-n}(x,\xi ,-1)\,\d\xi 
=-\tfrac 1m \int_{|\xi
|=1}b_{-n}(x,\xi )\, \d S(\xi ).
\tag1.26
$$
\comment
In fact, the identification of the contributions from
$q^h_{-m-n}(x,\xi ,-1)$ and $-\frac1m b_{-n}(x,\xi )$ holds {\it on each
ray} $\{t\eta  \mid t\ge 0\}$, $\eta  \in S^{n-1}$, cf.\ {\rm (1.25)}
(holds microlocally in this sense).
\endcomment
\endproclaim 
 \example{Remark 1.5} In this application of Lemma 1.3,
$f(\xi ,t)=q^h_{-m-n}(x,\xi ,t)$ is not only integrable at $\xi =0$ but
continuous there, for $t\ne 0$. Then for any $\xi \in S^{n-1}$, 
$f(\xi ,t)=|t|^{-1-n/m}f( |t|^{-1/m}\xi , -1)$ where $f(
|t|^{-1/m}\xi , -1)\to f(0,-1)$ for $t\to-\infty $, assuring that the
integrals in (1.24) exist. We can then say that 
the identification of the contributions from
$q^h_{-m-n}(x,\xi ,-1)$ and $-\frac1m b_{-n}(x,\xi )$ holds {\it on each
ray} $\{s\xi  \mid s\ge 0\}$, $\xi \in S^{n-1}$
(holds microlocally in this sense).
\endexample

\subhead 2. The trace defect formulas for closed manifolds \endsubhead

Let $A$ be a classical pseudodifferential of order $\sigma \in\Bbb R$,
and let $P$ be as in the preceding section; we now assume
for convenience that
$
m>n+\sigma $. 

It was shown
in \cite{GS1, Th\. 2.7} ($m$ integer $>0$) and 
\cite{Lo}, 
\cite{GH} ($m\in\rp$), that the kernel of $A(P-\lambda )^{-1}$
calculated in local coordinates has an
expansion on the diagonal, implying a trace expansion by integration
of the fiber trace in $x$:
$$
\aligned
K( A(P-\lambda )^{-1},x,x)
&\sim
\sum_{ j\ge 0  }   c_{ j}(x)(-\lambda ) ^{\frac{\sigma +n -j}m-1}+ 
\sum_{k\ge 0}\bigl(   c'_{ k}(x)\log (-\lambda ) +  c''_{
k}(x)\bigr)(-\lambda ) ^{{ -k}-1},\\
\Tr\bigl( A(P-\lambda )^{-1}\bigr)
&\sim
\sum_{ j\ge 0  }   c_{ j}(-\lambda ) ^{\frac{\sigma +n -j}m-1}+ 
\sum_{k\ge 0}\bigl(   c'_{ k}\log (-\lambda ) +  c''_{
k}\bigr)(-\lambda ) ^{{ -k}-1}.
\endaligned \tag2.1$$
Here $\lambda \to\infty $ on rays in 
an open subsector $V$ of $\Bbb
C$ containing $\Bbb R_-$.
It is convenient to assume that the operators are represented, via
local coordinate systems, as a finite sum of pieces acting separately
in a system of disjoint open sets in $\Bbb R^n$ (as e.g.\ in
\cite{G2, Sect.\ 1}), so that
we get the 
trace simply by integrating over $\Bbb R^n$.

The $ c'_k(x)$ vanish when $\sigma \notin \Bbb Z$.
We shall {\it define}
$$  c_{n+\sigma }(x)=0,\quad c_{n+\sigma }=0,\quad\text{ if }n+\sigma
\notin \Bbb N;\tag2.2$$ 
 then $  c_{n+\sigma }(x)$ and $c_{n+\sigma }$ have a meaning for any
$\sigma $.  (We denote $\{0,1,2,\dots\}=\Bbb N$.)

The coefficient of $(-\lambda )^{-1}$ in (2.3) will be denoted
$C_0(A,P)$;
$$
C_0(A,P)=   c_{\sigma +n}+ 
c''_0.
\tag2.3
$$
Corresponding to (2.1), the generalized zeta
function $\zeta (A,P,s)$, defined as 
$\Tr(AP^{-s})$ for large $\operatorname{Re}s$, has a meromorphic
extension to $\Bbb C$ with poles at the points $(j-n)/m$, with
Laurent coefficients directly related to the coefficients in the
expansion (2.1). In particular, $C_0(A,P)$ equals the
coefficient of $s^0$ plus $\Tr(A\Pi _0)$, cf.\ (1.4)ff.

It is well-known that $C_0(A,P)$ is in general nonlocal in the sense
that it
depends on the full structure of $A$, not just its homogeneous
symbols. However, when $A'$ and $P'$ are another pair of similar
operators, one can
show that 
$$
C_0(A,P)-C_0(A,P') \text{ and }C_0([A,A'],P)\text{ are local}\tag
2.4$$
(depend on a finite set of
strictly homogeneous symbol terms of $A$, $A'$, $P$ and $P'$);  
in this sense, $C_0(A,P)$ is a {\it
quasi-trace} on the classical $\psi $do's $A$. $C_0(A,P)$
is called a
regularized trace or weighted trace by other authors.
Explicit formulas for the trace defects in (2.4) were shown by
Okikiolu \cite{O}, 
Kontsevich and Vishik \cite{KV}, and Melrose and Nistor \cite{MN}:
$$\gather
C_0(A,P)-C_0(A,P')=-\tfrac 1m \operatorname{res}(A(\log P-\log P')),\tag2.5
\\
C_0([A,A'],P)=-\tfrac 1m \operatorname{res}(A[A',\log P]).\tag2.6
\endgather$$
Here Okikiolu proved (2.5) by an exact symbol calculation passing via
the symbols of the complex powers $P^{-s}$ and $(P')^{-s}$,
and Kontsevich and Vishik proved it by use of their calculus
of weakly holomorphic $\psi $do families. Melrose and Nistor
showed both (2.5) and (2.6) on the basis of the theorem of Guillemin
on holomorphic families \cite{Gu} (we have reconstructed a proof
based on this idea in \cite{G2, pf.\ of Prop.\
3.1}).
In all these cases, the logarithm $\log P$ comes up as a result of a
differentiation of $P^{-s}$ with respect to $s$.

Our present aim is to show how the formulas (2.5)--(2.6) can be found
directly from the knowledge of the resolvent expression $A(P-\lambda
)^{-1}$, without worrying about the construction of $P^{-s}$. (This
is important for generalizations to other types of manifolds.) We show
that in fact the full operator $\log P$ plays a very minor role; its
symbol comes in only 
because of the jump across the negative real axis as in Lemma 1.2.

\medskip

Let $P$ and $P'$ be auxiliary operators of order $m$ with resolvents
$Q_\lambda =(P-\lambda )^{-1}$, $Q'_\lambda =(P'-\lambda )^{-1}$
(symbols $q$ resp.\ $q'$), and 
consider the symbol $s(x,\xi ,\lambda )$ of $$
S_\lambda =A(Q_\lambda -Q '_\lambda )\tag2.7
$$ in local
coordinates. Much as in Lemma 1.1, we can show:

\proclaim{Proposition 2.1} The symbol $s$ of $S_\lambda =A(Q_\lambda
-Q'_\lambda )$ has an expansion in
strictly homogeneous terms plus a remainder:$$
s(x,\xi ,\lambda )=\sum_{0\le j<\sigma +m+n}s^h_{-m-j}(x,\xi ,\lambda
)+s'_{-2m-n}(x,\xi ,\lambda ),\tag2.8 
$$
where the $s^h_{\sigma -m-j}$ and $s'_{-2m-n}$ are integrable in $\xi
$ for $\lambda \ne 0$, and
$\int s'_{-2m-n}\,\d\xi $ is $O(|\lambda |^{-2+\varepsilon })$, any
$\varepsilon >0$.
 Consequently, $K(S_\lambda ,x,x) $ and the trace $\Tr S_\lambda $ have the expansions$$
\aligned
K(S_\lambda ,x,x)
&=
\sum_{ j<\sigma + m+n  }  \tilde s_{ j}(x)(-\lambda ) ^{\frac{n
+\sigma -j}m  -1}+
O(|\lambda |^{-2+\varepsilon }),\\
\Tr S_\lambda 
&=
\sum_{ j<\sigma +m+n }  \tilde s_{ j}(-\lambda ) ^{\frac{n +\sigma -j}m  -1}+
O(|\lambda |^{-2+\varepsilon }),\text{ where }\\
  \tilde s_j(x)&=\int_{\Bbb R^n}s^h_{-m-j}(x,\xi ,-1)\,\d\xi ,\quad 
  \tilde s_j=\int \tr \tilde s_j(x)\,dx ,\quad 
\text{ for } j<\sigma +m+n.
\endaligned
\tag 2.9
$$

In particular, when $n+\sigma \notin \Bbb N$, there is no term
with $(-\lambda )^{-1}$ in the expansion of  $\Tr S_\lambda $, and $$
C_0(A,P)-C_0(A,P')=0.\tag2.10
$$ 
When $n+\sigma \in\Bbb N$, the coefficient of $(-\lambda )^{-1}$ in
$\Tr S_\lambda $ equals  
$$
C_0(A,P)-C_0(A,P')=\int \tr \tilde s_{n+\sigma }(x)\,dx.\tag2.11
$$ 
\endproclaim 

\demo{Proof}
We use again the analysis of the resolvent symbol
recalled in Section 1. The composition with $A$ in front leads to
terms of the form (1.5) where the $\lambda $-independent coefficients
now furthermore contain information from the symbol $a$ of $A$. 
Consider$$
q-q'\sim
q_{-m}(x,\xi ,\lambda )-q'_{-m}(x,\xi ,\lambda )+
 \sum_{j\ge 1}(q_{-m-j}(x,\xi ,\lambda
-q'_{-m-j}(x,\xi ,\lambda )).
$$
All the terms in the sum over $j\ge 1$ are finite sums of
expressions as in (1.5),
containing at least two principal resolvent factors $q_{-m}$, resp\.
$q'_{-m}$. Moreover,
$$
\multline q_{-m}-q'_{-m}=(p_m-\lambda )^{-1}-(p'_m-\lambda )^{-1}\\=
(p_m-\lambda )^{-1}(p'_m-p_m)(p'_m-\lambda )^{-1}
=q_{-m}(p'_m-p_m) q'_{-m},\endmultline
$$
showing that it also contains two principal resolvent factors
($q_{-m}$ and $q'_{-m}$) together with a $\lambda $-independent 
factor.
Then an application of the standard composition rule 
gives that the homogeneous terms in the
symbol $s=a\circ(q-q')$ of $S_\lambda $
are finite sums of expressions that are a
slightly generalized version of (1.5) 
where some of the factors $q_{-m}$ may be replaced by $q'_{-m}$. The
important observation is that there are at least two such factors in
each term. Then, taking the order and homogeneity degrees into
account, we see that
$
s(x,\xi ,\lambda )\sim \sum_{j\ge 0}s_{\sigma -m-j}(x,\xi
,\lambda )$ satisfies $$
\aligned
|s^h_{\sigma -m-j}| &\le c(|\xi |^m+|\lambda |)^{-2}|\xi
|^{\sigma +m-j},\\
|s^h_{\sigma -m-j}-s_{\sigma -m-j}| &\le c|\lambda |^{-2}(1+|\xi
|^{\sigma +m-j}), \text{ supported in }|\xi |\le 1,\text{ any }j,\\
|s-\sum_{j<N}s_{\sigma -m-j} |&\le c(1+|\xi |^m+|\lambda |)^{-2}(1+|\xi
|)^{\sigma +m-N},\text{ any }N.
\endaligned
\tag2.12
$$
For $j<\sigma +m+n$, the two first expressions are
integrable in $\xi $. The remainder $s-\sum_{j<\sigma +m+n}s_{\sigma
-m-j}$ is seen as in the treatment of (1.11) to
be  $O((1+|\lambda |)^{-2+\varepsilon }(1+|\xi
|)^{-n-\delta  '})$ with $\delta '>0$ and $\varepsilon $ arbitrarily
small, here $\varepsilon $ can be taken $=0$ if $\sigma +m\notin\Bbb Z$.  
 This shows the first part of the lemma, and the second part follows by
integration, first in $\xi $ and then (for the fiber trace) in $x$.

For the third part,
observe that there is no term $c(-\lambda )^{-1}$  in (2.9) when
$n+\sigma\notin \Bbb N$. When $n+\sigma\in \Bbb N$,
the coefficient 
of  $(-\lambda )^{-1}$ is found from (2.9) for $j=n+\sigma $.\qed
\enddemo 

Note that all the
indicated coefficients are local, and that there is no $(-\lambda
)^{-1}\log(-\lambda )$ term as in (2.3).

We can now show (2.5) in a precise form, by a calculation as in
Section 1.
For this we consider 
$$
F=A(\log P-\log P').
$$ 
Since the logarithmic
terms in the symbols of $\log P$ and $\log P'$ cancel out (cf.\
(1.14)), it is a classical $\psi $do of order $\sigma $; we 
denote its symbol by $f(x,\xi )$. 
When we define $F$ by the
formula
$$
\aligned
F&=A(\log P-\log P')=A\lim_{s\to 0}\tfrac i{2\pi }\int_{\Cal C}\lambda
^{-s}\log\lambda \,\bigl((P-\lambda )^{-1}-(P'-\lambda
)^{-1}\bigr)\,d\lambda \\
&=\tfrac i{2\pi }\int_{\Cal C}\log\lambda \,S_\lambda \,d\lambda , 
\endaligned\tag2.13
$$
then in local coordinates, its symbol
is found termwise from the symbol of
$S_\lambda $ by the formulas
$$
f_{\sigma -j}(x,\xi )=\tfrac i{2\pi }\int_{\Cal C'}\log\lambda
\,s_{\sigma -m-j}(x,\xi ,\lambda )\, d\lambda ,\tag2.14
$$
where $\Cal C'$ is a closed curve in $\Bbb C\setminus\crm$ encircling
the eigenvalues of 
$p_m$ and $p'_m$. This follows from the calculations of the terms in
$\log P$ and $\log P'$ described e.g.\ in \cite{O}, and the composition 
rule for $\psi $do's.

When $n+\sigma \notin \Bbb N$, there is no term of degree $-n$, so the
noncommutative residue of $F$ is zero. When $n+\sigma \in\Bbb N$, it
is determined by
$$
\operatorname{res}F=\int_{X}\int_{|\xi |=1}\tr f_{-n}(x,\xi
)\,\d S(\xi )dx.\tag2.15
$$

\proclaim{Theorem 2.2} Let $P$ and $P'$ be classical $\psi $do's of
order $m>0$ and such that the principal symbol has no eigenvalues on
$\crm$, let $A$ be a 
classical $\psi $do of order $\sigma $, and let
\linebreak$S_\lambda =A((P-\lambda )^{-1} 
-(P'-\lambda )^{-1})$ and 
$F=A(\log P-\log P')$ with symbols $s$ resp.\ $f$. Assume that
$m>n+\sigma $.

Consider the case $n+\sigma \in\Bbb N$.
The formula {\rm (2.5)} is valid, and it holds pointwise, in
that
$$
C_0(A,P)-C_0(A,P')=\int_X \tr \tilde s_{n+\sigma }(x)\,dx
=-\tfrac1m\operatorname{res}(A(\log P-\log P'))
\tag2.16$$ 
where, for each $x$, in local coordinates,
$$
\tilde s_{n+\sigma }(x)=\int_{\Bbb R^n}s^h_{-m-n}(x,\xi ,-1)\,\d\xi 
=-\tfrac 1m \int_{|\xi
|=1}f_{-n}(x,\xi )\, \d S(\xi ).
\tag2.17
$$
\comment
In fact, the identification of the contributions from
$s^h_{-m-n}(x,\xi ,-1)$ and $-\frac1m f_{-n}(x,\xi )$ holds on each
ray $\{t\eta  \mid t\ge 0\}$, $\eta \in S^{n-1}$
(holds microlocally in this sense).
\endcomment

When $n+\sigma \notin \Bbb N$, the identities hold trivially (with zero
values everywhere).

\endproclaim 

\demo{Proof} The proof consists of rewriting $\int_{\Bbb
R^n}s^h_{-m-n}(x,\xi ,-1)\,\d\xi $ in the same way as we did with the
integral of $q^h_{-m-n}$ in Section 1:
$$\aligned
\int_{\Bbb R^n}s^h_{-m-n}(x,\xi ,-1)\,\d\xi 
& =-\tfrac1m \int_{|\xi |=1}\tfrac i{2\pi }\int_{\Cal C'}\log \lambda
\,s^h_{-m-n}(x,\xi ,\lambda )\,d\lambda 
\d S(\xi )\\
&=-\tfrac 1m \int_{|\xi
|=1}f_{-n}(x,\xi )\, \d S(\xi ),
\endaligned\tag2.18$$
where the first equation follows from Lemmas 1.2 and 1.3, and the
second equation follows from (2.14).\qed
\enddemo

There is a related proof of the other trace defect formula, (2.6).
We here consider $A$ of order $\sigma $, $A '$ of order $\sigma '$
and $P$ as before, now assuming for convenience that
$
\sigma +\sigma '+m>n$. 

Here we first observe that by cyclic permutation,$$
\aligned
\Tr([A,A']Q_\lambda )&=\Tr( AA'Q_\lambda )-\Tr(AQ_\lambda
A')=\Tr(A[A',Q_\lambda]),\text{ where}\\
A[A',Q_\lambda ]&=A(Q_\lambda (P-\lambda )A'Q_\lambda -Q_\lambda A'(P-\lambda
)Q_\lambda )=AQ_\lambda [P,A']Q_\lambda .
\endaligned
\tag2.19
$$
Let 
$$
\align 
T_\lambda &=[A,A']Q_\lambda ,\tag2.20\\
 R_\lambda &=A[A',Q_\lambda ]=AQ_\lambda [P,A']Q_\lambda .\tag2.21
\endalign
$$
The traces of $T_\lambda $ and $R_\lambda $ are identical, and the operators both have
order $\sigma +\sigma 
'-m$. It is seen from the second formula for $R_\lambda $ that the homogeneous
terms $r_{\sigma +\sigma '-m-j}$ in its symbol 
$r$ are finite sums of terms of the form (1.5) with at least two
factors $q_{-m}$, so that 
the strictly homogeneous symbols $r^h_{\sigma +\sigma '-m-j}$ are
integrable in $\xi $ at $\xi =0$ for $j<\sigma +\sigma '+m+n$.

We then find very similarly to the study of $S_\lambda $ that the diagonal kernel
of $R_\lambda $ has an expansion 
$$
K(R _\lambda ,x,x)
=
\sum_{ j<\sigma +\sigma '+m+n }   \tilde r_{ j}(x)(-\lambda )
^{\frac{n+\sigma +\sigma ' -j}m  -1}+ 
O(|\lambda |^{-2+\varepsilon }),\tag2.22
$$
where 
$$
\tilde r_j(x)=  \int_{\Bbb R^n}r^h_{\sigma +\sigma '-m-j}(x,\xi ,-1)\,\d\xi .\tag2.23
$$
When $n+\sigma+\sigma '\notin \Bbb N$, there is no term $c(-\lambda
)^{-1}$  in (2.22), hence no such term in the trace expansion of
$R_\lambda $.
Since this is the same as that of $T_\lambda $, the term is also missing from
$\Tr T_\lambda $, so $C_0([A,A'],P)=0$. When $n+\sigma+\sigma '\in \Bbb N$,
the coefficient
of  $(-\lambda )^{-1}$ in (2.22) equals (2.23) with $j=n+\sigma
+\sigma '$, i.e., 
$$
\tilde r_{n+\sigma +\sigma '}(x)=  \int_{\Bbb R^n}r^h_{ -m-n}(x,\xi ,-1)\,\d\xi .\tag2.24 
$$
Then the coefficient of $(-\lambda )^{-1}$ in the expansion of $\Tr
R_\lambda $ equals the integral in $x$ of the fiber trace of this (collecting
the contributions 
from local coordinate systems), and since $\Tr T_\lambda $ has the same
expansion, we can conclude that
$$
C_0([A,A'],P)=\int_X \tr \tilde r_{n+\sigma +\sigma '}(x)\,dx.\tag2.25
$$ 

On the other hand, we consider $H=A[A',\log P]$, observing that it is
a classical $\psi $do of order $\sigma +\sigma '$ in view of (1.14).
Here,
$$
\aligned
H&=A(A'\log P-\log P\, A')=A\lim_{s\to 0}\tfrac i{2\pi }\int_{\Cal C}\lambda
^{-s}\log\lambda \,(A'Q_\lambda -Q_\lambda A')\,d\lambda\\
&=\tfrac i{2\pi }\int_{\Cal C}
\log\lambda \,R_\lambda \,d\lambda.
\endaligned
\tag2.26
$$
The symbol $h(x,\xi )$ is found termwise in local coordinates from
the symbol of $R_\lambda $ by the formulas
$$
h_{\sigma +\sigma '-j}(x,\xi )=\tfrac i{2\pi }\int_{\Cal C'}\log\lambda
\,r_{\sigma +\sigma '-m-j}(x,\xi ,\lambda )\, d\lambda ,\tag2.27
$$
where $\Cal C'$ is a closed curve in $\Bbb C\setminus\crm$ encircling
the eigenvalues of 
$p_m$.

When $n+\sigma +\sigma '\notin \Bbb N$, there is no term of degree $-n$ so the
noncommutative residue of $H$ is zero. When $n+\sigma +\sigma '\in\Bbb N$, it
is determined by
$$
\operatorname{res}H=\int_{X}\int_{|\xi |=1}\tr h_{-n}(x,\xi
)\,\d S(\xi )dx.\tag2.28
$$

We then get:

\proclaim{Theorem 2.3} With $P$ and $A$ as in Theorem {\rm 2.2}, let
$A'$ be a  
classical $\psi $do of order $\sigma '$, and let $R_\lambda
=A[A',(P-\lambda )^{-1}])$ and 
$H=A[A',\log P]$ with symbols $r$ resp.\ $h$. Assume that $m>n+\sigma
+\sigma '$.

Let $n+\sigma +\sigma '\in\Bbb N$.
The formula {\rm (2.6)} is valid, and it holds pointwise, in
that
$$
C_0([A,A'],P)=\int_X \tr \tilde r_{n+\sigma +\sigma '}(x)\,dx
=-\tfrac1m\operatorname{res}(A[A',\log P])
\tag2.29$$ 
where, for each $x$, in local coordinates,
$$
\tilde r_{n+\sigma +\sigma '}(x)=\int_{\Bbb R^n}r^h_{-m-n}(x,\xi ,-1)\,\d\xi 
=-\tfrac 1m \int_{|\xi
|=1}h_{-n}(x,\xi )\, \d S(\xi ).
\tag2.30
$$

When $n+\sigma +\sigma '\notin \Bbb N$, the identities hold trivially
(with zero 
values everywhere).

\endproclaim 

\demo{Proof} 
The
identity follows from (2.27) together with Lemmas 1.2 and 1.3, in the
same way as in Theorem 2.2.\qed 
\enddemo 

\example{Remark 2.4} The observation in Remark 1.5 on the microlocal
identification extends to the
formulas in Theorems 2.2 and 2.3.
\endexample

\comment
\example{Remark 2.5}
Also here, one could proceed instead by appealing to considerations on
regularity numbers as in \cite{G1}, using that $R$ according to the
last expression in (2.21) has regularity number $\nu =\sigma +\sigma
'+m$
(by a generalization of \cite{G1, Lemma 2.1.6}).
\endexample
\endcomment

\subhead  3. The first trace defect formula for manifolds with boundary
\endsubhead
  
We shall now discuss extensions of the above
results to pseudodifferential boundary operatos ($\psi $dbo's) of
Boutet de Monvel's type in the case of manifolds with boundary.

Consider a compact $n$-dimensional $C^\infty $ manifold $X$ with
boundary $\partial X=X'$, and a hermitian $C^\infty $ vector bundle
$E$ over $X$.
Let $A=P_++G$ be an operator of order $\sigma $ belonging to the
calculus of Boutet de Monvel \cite{B}, acting on sections of $E$.
Here $P$ is a classical $\psi $do satisfying the transmission
condition at $\partial X$ and $G$ is a singular Green operator
(s.g.o.)\ of
class 0 with polyhomogeneous symbol. (More details can be found e.g.\
in \cite{B, G1}). When $P\ne 0$, we must assume
$\sigma \in\Bbb Z$ because of the requirements of the transmission
condition; when $P=0$ it is straightforward to allow $\sigma \in\Bbb R$. 
For the results in Section 4, $P$ is moreover assumed to be of normal
order $\le 0$ (its symbol is bounded in $\xi _n$, the boundary
conormal variable).

As auxiliary operator we take an elliptic differential operator $P_1$
of order $m>0$ with scalar principal symbol $p_{1,m}$ taking no
values on $\rmi$; 
so $(p_{1,m}(x,\xi )-\lambda )^{-1}$ is defined for $\lambda $ in a
sector $V$ around 
$\rmi$, for all $x$, all $\xi $ with $|\xi |+|\lambda |\ne 0$.
$P_1$ can be assumed to be given on a larger
boundaryless $n$-dimensional compact manifold $\widetilde X$ in which $X$ is
smoothly imbedded, acting in a bundle $\widetilde E$ extending $E$
and with the same ellipticity properties there.
We set$$
Q_\lambda =(P_1-\lambda )^{-1}\tag 3.1
$$
on $\widetilde X$; it is defined except for a discrete subset of
$\Bbb C$; in particular it exists for large $\lambda $ in the sector $V$.

For the case where $m=2$ and $P_1$ is strongly elliptic, defining the
Dirichlet realization $P_{1,\operatorname{D}}$, we showed  
in a joint work with Schrohe  \cite{GSc1}  that there is a resolvent
trace expansion when $N>(\sigma +n)/2$: 
$$
\Tr(A(P_{1,\operatorname{D}}-\lambda )^{-N})\sim \sum_{j\ge
0}\!\tilde c_{j}(-\lambda)
^{\frac{n+\nu - j}{2}-N}
+\sum_{k\ge
0}(\tilde c'_k\log (-\lambda )+\tilde c''_k)(-\lambda )
^{-\frac k2-N}, \tag3.2
$$   
valid for $\lambda \to \infty $ in  $V$. It was used there
to show that the coefficient $\tilde c'_0$ is proportional
to the
noncommutative residue of $A$, as introduced by Fedosov, Golse,
Leichtnam and Schrohe in \cite{FGLS}.

The proofs in \cite{GSc1} were formulated only for $\sigma \in\Bbb
Z$; but for more general $\sigma \in\Bbb R$, they carry over without
difficulty to the case $A=G$. In particular, if $\sigma \in\Bbb
R\setminus\Bbb Z$, the coefficients $\tilde c'_k$ vanish (since the
$\psi $do's on the boundary obtained by reduction of
$A(P_{1,\operatorname{D}}-\lambda )^{-N}$ are polyhomogeneous of noninteger
order). The identification of $\tilde c'_0$ with a noncommutative residue
then holds with$$
\operatorname{res}(G)=0, \text{ when }\sigma \notin\Bbb Z.\tag3.3
$$

As usual we define $$
C_0(A,P_{1,\operatorname{D}})= \tilde c_{n+\sigma }+\tilde c''_0, 
$$
where $\tilde c_{n+\sigma }$ is defined to be 0 if $n+\sigma
\notin\Bbb N$.
By a precise analysis of the terms entering in trace
expansions like (3.2), we showed in \cite{GSc2}
that the
functional $C_0(A,P_{1,\operatorname{D}})$ has quasi-trace
properties as in (2.4); moreover, we singled out some cases where it has a
value independent of $P_1$ and vanishes on commutators, so that it can
be regarded as a canonical trace in a similar sense as that of \cite{KV}.

It is shown in \cite{GSc2} that the singular Green part $G_\lambda ^{(N)}$ of
$(P_{1,\operatorname{D}}-\lambda )^{-N}=(Q_\lambda ^{N})_++G_\lambda ^{(N)}$
contributes only locally to  
$C_0(A,P_{1,\operatorname{D}})$. It has an interest to consider the
composition $A(Q^N_\lambda )_+$ alone; it likewise has an expansion 
$$
\Tr(A(Q_\lambda ^N)_+)\sim  \sum_{j\ge
0}\!\tilde a_{j}(-\lambda)
^{\frac{n+\nu - j}{2}-N}
+\sum_{k\ge
0}(\tilde a'_k\log (-\lambda )+\tilde a''_k)(-\lambda )
^{-\frac k2-N}, \tag3.4
$$   
where $\tilde a'_0=\frac 1m \operatorname{res}(A)$, and the
coefficient of $(-\lambda )^{-N}$,
$$
C_0(A,P_{1,+})=\tilde a_{n+\sigma }+\tilde a''_0 \tag3.5
$$
is a quasi-trace on the $\psi $dbo's (by the results of \cite{GSc2}).

One may remark that in an associated zeta function formulation, the
consideration of $(Q^N_\lambda )_+$ alone 
corresponds to considering compositions with $(P_1^{-s})_+$ alone, 
where $(P_1^{-s})_+$ is another
family of operators than $(P_{1,\operatorname{D}})^{-s}$; both
families have the property that they equal $I$ when $s=0$. 

But actually these complex powers lie outside the
Boutet de Monvel calculus (when $s\notin\Bbb Z$). There is a
description in \cite{G1, Sect.\ 4.4} of negative powers
($\operatorname{Re}s>0$), showing how the s.g.o.\ part satisfies some
but not all the standard estimates. But they have not, to
our knowledge, been 
successfully described a holomorphic family in some sense where
results like that of Guillemin \cite{Gu, Th.\ 7.1} for closed
manifolds could be applied to generalize the trace defect formulas
(2.5)--(2.6). (The use of Guillemin's result is explained e.g.\ in
\cite{G2, pf.\ of Prop.\ 3.1}.)

Even if one avoids dealing with complex powers, there is still the
problem in generalizing the formulas (2.5)--(2.6) that
logarithms of $\psi $dbo's have not been studied, and do
not in general belong to the Boutet de Monvel calculus. However,
$(\log P_1-\log P_2)_+$ does belong there when $P_1$ and $P_2$ are
two choices of the auxiliary elliptic operator (of order $m$), thanks
to the cancellation of 
logarithms resulting from  (1.14). But $[A',(\log P_1)_+]$ does not
so, except in trivial cases. 

We shall show a generalization of (2.5) in this section, and treat
(2.6) in the following section.

\medskip
The papers \cite{GSc1, GSc2} used the refined calculus of Grubb and
Seeley \cite{GS1},
which allows obtaining complete trace expansions (with remainders $O(\lambda
^{-M})$, any $M$).

Presently we shall use the cruder (but more generally applicable)
calculus from the book \cite{G1} to achieve our result, building also
on the insight gained in Sections 1 and 2. Notably, we are avoiding
some technical challenges by restricting the attention to the trace of
$AQ_{\lambda ,+}$, without an s.g.o.\ term $AG_\lambda $ coming from
a boundary condition on $P_1$.

An advantage is that we can allow rather general auxiliary operators
$P_1$ of higher order, with no conditions on root multiplicities in the
principal symbol. (In \cite{GSc1, GSc2}, the order 2 assured
well separated roots in $\xi _n$, one in
each complex half-plane.) On the other hand, the theory we presently use
gives trace expansions with a finite number of terms only (plus a
remainder); but this turns out to be just sufficient for studying the
trace defect formulas. 

\medskip
Let us first recall some elements of the theory. As usual, $\ang\xi $
stands for $(1+|\xi |^2)^{\frac12}$; moreover, it is convenient to
denote $\ang{(\xi ,\mu )}=\ang{\xi ,\mu }$ and use the sign $\leg$
as shorthand for ``$\le$ a constant times''.

A  $\psi $do symbol $s(x,\xi ,\mu )$ on $\Bbb R^n$ depending on the
parameter $\mu 
\in\crp$ is said to be {\it of order $d $ and regularity $\nu $} ($d,\nu
\in \Bbb R$) with uniform estimates (\cite{G1, Def.\ 2.1.1}), when it
satifies, for all indices $\alpha ,\beta ,j$: 
 $$
|D_x^\beta D_\xi^\alpha D_\mu^js(x,\xi,\mu)| \le 
c\bigl(\ang\xi^{\nu-|\alpha|}+\ang{\xi,\mu}^{\nu-|\alpha|}\bigr)
\ang{\xi,\mu}^{d-\nu-j}
\tag3.6
 $$
 for $(x,\xi,\mu)\in \Bbb R^{2n}\stimes \crp$, with constants $c$ depending on
the indices. It is then said to be {\it polyhomogeneous}, when it
furthermore has an expansion $s(x,\xi,\mu)\sim\sum_{l\in\Bbb
N}s_{d-l}(x,\xi,\mu) 
$ in terms $s_{d-l}$ that are homogeneous in $(\xi ,\mu )$ of degree
$d-l$ for $|\xi |\ge 1$, such that $s-\sum_{l<M}s_{d-l}$ is of order
$d-M$ and regularity $\nu -M$, for all $M\in\Bbb N$. Note that in
(3.6), $\ang\xi ^{\nu -|\alpha |}$ can be left out when $|\alpha |\le
\nu $, and $\ang{\xi ,\mu } ^{\nu -|\alpha |}$ can be left out when
$|\alpha |\ge \nu $.

For such symbols we have:

\proclaim{Lemma 3.1} Let $s(x,\xi ,\mu )$ be polyhomogeneous of order
$d$ and regularity $\nu $, $d$ and $\nu \in\Bbb R$. Write $\nu =\nu '+\delta $ with $\nu
'$ integer and $\delta \in \,]0,1]$. Then 
$$\aligned
&|s_{d-l}(x,\xi,\mu)|\leg\ang{\xi,\mu}^{d-l}\text{ for
}l\leq \nu ,\text{ with}\\
&|s_{d-l}^h(x,\xi,\mu)|\leg |(\xi,\mu)|^{d-l},\\
&|s_{d-l}-s_{d-l}^h|\leg\ang{\xi,\mu}^{d-\nu }\text{ for
}|(\xi,\mu)|\geq c>0;\endaligned\tag3.7
 $$
 and the next terms with $l\le \nu +n$ are estimated by
 $$\aligned
&|s_{d-l}(x,\xi,\mu)|\leg 
\ang\xi^{\nu -l}\ang{\xi,\mu}^{d-\nu },\text{ for }\nu <l<\nu +n,\text{
with}\\
&|s_{d-l}^h(x,\xi,\mu)|\leg
|\xi|^{\nu -l}|(\xi,\mu)|^{d-\nu };\endaligned\tag3.8
$$
so that altogether
 $$s(x,\xi,\mu)=s_{d}^h(x,\xi,\mu)+\dots+s_{d-\nu '-n+1}^h
(x,\xi,\mu)+s'(x,\xi,\mu)\tag3.9
 $$
 where $s'=s''+s'''$, satisfying for $\mu\geq c_0>0$,
 $$\aligned
&|s''|\leg|\xi|^{\delta -n}|\mu|^{d-\nu }\chi(\xi),\text{ with }\chi\in
C_0^\infty(\Bbb R^n),\chi(\xi)=1\text{ for }|\xi|\leq1,\\
&|s'''|\leg\ang\xi^{\delta -n-1}\ang{\xi,\mu}^{d-\nu }.\endaligned\tag3.10
 $$ 
\endproclaim 

\demo{Proof} The proof, given for particular choices of $d$ and $\nu $
in \cite{G1, (3.3.35)ff.\ and (3.3.69)ff.}, extends to the general
situation: The first two lines in 
(3.7) follow readily from the definitions. The third line is
less obvious; it is shown in \cite{G1, Lemma 2.1.9 2$^\circ$} (by integration
of an estimate of a high enough derivative). (3.8) follows easily from the
definitions (one may consult \cite{G1, Lemma 2.1.9 1$^\circ$}). Then
(3.9) follows in view of (3.10), where $s''$ collects the differences
between homogeneous and 
strictly homogeneous symbols and $s'''$ is the remainder
$s-\sum_{l\le \nu '+n}s_{d-l}$.\qed 
\enddemo 

As in \cite{G1, Th.\ 3.3.5 and 3.3.10} we can use the lemma to get a
diagonal kernel expansion with $n+\nu '$ precise terms:

\proclaim{Lemma 3.2} When $d<-n$ in Lemma {\rm 3.1}, the kernel
$K(S_\mu ,x,y )$ of $S_\mu =\op (s(x,\xi ,\mu ))$ is continuous and
has an expansion on 
the diagonal: 
$$
K(S_\mu ,x,x )=\sum_{0\le l<n+\nu }\tilde s_{l}(x)\mu
^{n+d-l}+\tilde s'(x,\mu );\tag3.11
$$
here 
$$
\tilde s_{l}(x)=
\int_{\Bbb R^n}s^h_{d-l}(x,\xi ,1)\,\d\xi, 
\tag3.12
$$ 
and
 $\tilde s'(x,\mu )$ is $O(\mu ^{d-\nu +\varepsilon })$ for $\mu
\to\infty $, any $\varepsilon >0$. Here if $\nu \notin \Bbb Z$,
$\varepsilon $ can be left out.
\endproclaim 

\demo{Proof} This follows by integration of (3.9) in $\xi $.
For the terms $ s^h_{d-l}$ we use the homogeneity, replacing $\xi $
by $\eta =\mu ^{-1}\xi $:$$
\int_{\Bbb R^n}s^h_{d-l}(x,\xi ,\mu )\,\d\xi =\mu ^{d-l+n}\int_{\Bbb
R^n}s^h_{d-l}(x,\eta  ,1 )\,\d\eta , 
$$ 
and for $s'$ we use the estimates (3.10) 
(cf.\ \cite{G1, Lemma 3.3.6}).\qed
\enddemo 

\example{Remark 3.3}
The symbol spaces $S^{r,a}_{\operatorname{phg}}$ ($a\in\Bbb Z$)
defined in \cite{GS1} are somewhat more refined; they fit 
into the regularity classes as follows: Let $s(x,\xi ,\mu )$ belong
to $S^{r,a}_{\operatorname{phg}}\cap S^{r+a,0}_{\operatorname{phg}}$,
where $r+a< -n$. Then $s$ is of order $r+a$, and
$f(x,\xi ,\mu )=\mu ^{-a}s(x,\xi ,\mu )$ satisfies the requirements
(3.6), except those concerning $\mu $-derivatives, for being of order $r$ and
regularity $r$. The fact that $r+a<-n$ makes the symbol integrable in
$\xi $; the information on $f$ assures that the strictly homogeneous
terms $f^h_{r-l}$ are integrable at $\xi =0$ for $l<r+n$ and the
remainder $f-\sum_{l<n+r}f^h_{r-l}$ is integrable at $\xi =0$ (as
in Lemma 3.1). Then we get the diagonal expansion of the
kernel of $S_\mu =\op (s)$ as in Lemma 3.2:$$
K(S_\mu ,x,x )=\mu ^a\sum_{0\le l<n+r }\tilde s_{l}(x)\mu
^{n+r-l}+\tilde s'(x,\mu );\quad \tilde s'(x,\mu )=O(\mu
^{a+\varepsilon }),\tag3.13
$$
for $\mu \to\infty $, with locally determined coefficients $\tilde
s_l$. What the calculus of \cite{GS1} moreover gives for the symbols
in $S^{r,a}_{\operatorname{phg}}\cap S^{r+a,0}_{\operatorname{phg}}$
is a full expansion of the 
remainder:
$$
\tilde s'(x,\mu )\sim \mu ^a[\sum_{l\ge n+r}\tilde s_{l}(x)\mu
^{n+r-l}+\sum_{k\ge 0}(\tilde s_k'(x)\log\mu
+\tilde s''_k(x))\mu ^{-k}],
$$
with local coefficients $\tilde s_l(x)$, $\tilde s'_k(x)$ and global coefficients
$\tilde s''_k(x)$. Some of the $\tilde s''_k(x)$ may belong to the same
powers as 
coefficients $\tilde s_l(x)$, so the values of $a$ and $r$ are
important in the discussion of which terms are local.
\endexample

Besides $\psi $do's we must now deal with
singular Green operators.
Singular Green symbol-kernels $\tilde g(x',x_n,y_n,\xi ',\mu )$ of
{\it order} $d$ (degree $d-1$), regularity $\nu $ and
class $0$ satisfy estimates
$$
\|D_{x'}^\beta D_{\xi'}^\alpha D_\mu^j \tilde
g(x',x_n,y_n,\xi',\mu)\|_{L_2(\rp\stimes\rp)}
\leg(\ang{\xi '}^{\nu-|\alpha|}+\ang{\xi ',\mu }^{\nu -|\alpha |})
\ang{\xi ',\mu }^{d-\nu -j},\tag3.14 
$$
along with further estimates for $x_n^k
D_{x_n}^{k'}y_n^mD_{y_n}^{m'}\tilde g$ ($k,k',m,m'\in\Bbb N$) (cf.\
\cite{G1, Sect.\ 2.3}), and with a suitable definition of polyhomogeneity.

There is the following rule for normal traces of s.g.o.\ symbol-kernels: 

\proclaim{Lemma 3.4} When $\tilde g(x',x_n,y_n,\xi ',\mu )$ is a
singular Green symbol-kernel of order $d$, regularity $\nu $ and
class $0$, then the normal trace$$
s(x',\xi ',\mu )=\tr_n \tilde g=\int_0^\infty \tilde
g(x',x_n,x_n,\xi ',\mu )\, dx_n\tag3.15
$$
is a $\psi $do symbol on $\Bbb R^{n-1}$ of order $d$ and regularity
$\nu -\frac14$, polyhomogeneous if $\tilde g$ is so. 
\endproclaim 

\demo{Proof} This is shown in \cite{G1, pf.\ of Th.\ 3.3.9}
for $\nu \in \Bbb Z\cup (\Bbb Z+\frac12)$; $s$ is denoted $\ttilde g$
there. The first part of the proof
extends to all real $\nu \ge 1$ or $\le 0$; the loss of $\frac14$
stems from the negative cases (which occur in symbol terms of low
order). The last part of the proof, showing how $\nu =\frac12$ is
included by use of a derivative, extends to
general $\nu \in \,]0,1[$. (The considerations in \cite{G1} were
aimed at negative integer values of $d$, but all the arguments work
with arbitrary $d\in\Bbb R$ also.)\qed
\enddemo 

Combining Lemma 3.4 with Lemma 3.2 in dimension $n-1$, we find for
s.g.o.s of order $d<1-n$ and regularity $\nu $:

\proclaim{Lemma 3.5} Let $G_\mu $ be a $\mu $-dependent polyhomogeneous
singular Green operator on $\crnp$ of order $d<1-n$ and regularity $\nu
$. 
The normal trace $S_\mu =\tr_n G_\mu $ is a
$\psi $do on $\Bbb R^{n-1}$ of order $d$ and regularity $\nu -\frac14$,
whose kernel on the diagonal has an expansion in powers of $\mu $:
$$
K(S_\mu ,x',x')=\sum_{0\le l<n-1+\nu -\frac14}\tilde s_{l}(x')\mu
^{n-1+d-l}+O(\mu ^{d-\nu +\frac14 \,(+\varepsilon )});\tag3.16
$$
with $\varepsilon =0$ if $\nu -\frac14\notin\Bbb Z$, any small
$\varepsilon >0$ if $\nu -\frac14\in\Bbb Z$. Here
$$
\tilde s_{l}(x')=
\int_{\Bbb R^{n-1}}s^h_{d-l}(x',\xi ',1)\,\d\xi'. 
\tag3.17
$$ 

When the kernel has compact $(x',y')$-support, $G_\mu $ is trace-class and
the trace has an expansion with coefficients $\tilde s_l=\int\tr
\tilde s_l(x')\, dx'$:
$$
\Tr G_\mu =\sum_{0\le l<n-1+\nu -\frac14}\tilde s_{l}\mu
^{n-1+d-l}+O(\mu ^{d-\nu +\frac14 \;(+\varepsilon )});\tag3.18
$$
with $\varepsilon $ as above.
\endproclaim 

\demo{Proof} By Lemma 3.4, the operator family $S_\mu =\tr_nG_\mu $
satisfies the 
hypotheses of Lemma 3.2 in $n-1$ dimensions with $\nu $ replaced by
$\nu -\frac14$, this implies (3.16) with (3.17). Then (3.18) follows by
integration in $x'$.\qed
\enddemo 

Now let us turn to the specific operators we want to study.
Consider $A=P_++G$ of order $\sigma \in\Bbb R$ together with an auxiliary
elliptic operators $P_1$ of order 
$m>n+\sigma $. Recall that if $\sigma \in\Bbb R\setminus \Bbb Z$ we have
$P=0$.

The resolvent $Q_{1,\lambda }=(P_1-\lambda )^{-1}$
depends on $\lambda $ running in a sector $V$ around $\rmi$
in $\Bbb C$, where it is defined for large $\lambda $. We consider
$\lambda $ on each ray there, writing $-\lambda =\mu ^me^{i\omega }$, $\mu \ge
0$. Since $P_1$ is a differential operator, $Q_{1,\lambda }$ is of
regularity $+\infty $.  By \cite{G1, (2.1.13), (2.3.54)}, $A=P_++G$
enters in the 
parameter-dependent calculus as an operator of order and regularity
$\sigma $ (since $G$ is of class 0). Then the composed operator
$AQ_{1,\lambda ,+}$ is of order $\sigma -m$ and regularity $\sigma $,
in view of \cite{G1, Th.\ 2.7.7, Cor.\ 2.7.8} (no loss of
$\varepsilon $ regularity thanks to the mentioned theorem).

In the following, we work in a localized situation, as explained
e.g.\ in \cite{GSc2,  after (3.11)}.

For the $\psi $do 
$PQ_{1,\lambda }$ we already have a diagonal kernel expansion
(2.2) pointwise for $x\in\widetilde X$; integration of the fiber
trace over $X$ gives the trace expansion
 $$
\Tr((PQ_{1,\lambda })_+ )=
\sum_{ l\ge 0  }   c_{ l,+}(-\lambda ) ^{\frac{\sigma +n -l}m-1}+ 
\sum_{k\ge 0}\bigl(   c'_{ k,+}\log (-\lambda ) +  c''_{
k,+}\bigr)(-\lambda ) ^{{ -k}-1}.
\tag3.19
$$

Lemma 3.5 applied to the singular Green part $G_\lambda =AQ_{1,\lambda
,+}-(PQ_{1,\lambda })_+$ gives
the expansion
$$
\aligned
\Tr G_\lambda &=\sum_{0\le l<n-1+\sigma -\frac14 } b_{l}(-\lambda )
^{\frac{n-1+\sigma-l}m -1}+O(\lambda  ^{-1+\frac1{4m}(+\varepsilon )
})\\
&=\sum_{1\le j<n+\sigma -\frac14 } b'_{j}(-\lambda )
^{\frac{n+\sigma-j}m -1}+O(\lambda  ^{-1+\frac1{4m}\;(+\varepsilon )})
,\endaligned\tag3.20
$$
with $b'_j=b_{j-1}$; $\varepsilon $ equals $0$ when $\sigma
-\frac14\notin\Bbb Z$, 
and can be any small positive number when $\sigma -\frac14\in\Bbb Z$.
Here we first show the expansion on each ray, noting that for
$\lambda $ on the ray
 $\rmi$  we get (3.20); then the holomorphy
assures that the expansion is the same on the other rays
(as in \cite{GS1, Lemma 2.3}).
When $\sigma 
\notin\Bbb Z$, $G_\lambda =AQ_{1,\lambda ,+}$, so (3.20) shows its
trace expansion. When $\sigma \in\Bbb Z$,
addition of (3.19) and (3.20) gives: 
$$
\Tr (AQ_{1,\lambda ,+}) =\sum_{0\le l<n+\sigma  } c_{l}(-\lambda )
^{\frac{n+\sigma-l}m -1}+O(\lambda  ^{-1+\frac1{4m} }).\tag3.21
$$

This expansion does {\it not} show the appearance of a term
$c(-\lambda )^{-1}$. We shall obtain that by proving two things: 

1) When
$P_1$ is replaced by another auxiliary operator $P_2$ of order $m$,
then the difference of the traces $\Tr(AQ_{1,\lambda
,+})-\Tr(AQ_{2,\lambda ,+})$ has a better expansion:
$$
\Tr\bigl(A(Q_{1,\lambda }-Q_{2,\lambda })_+\bigr)
=\sum_{0\le j < n+\sigma +\frac14}d_{j}\,(-\lambda )
^{\frac{n+\sigma -j}m-1}+O(\lambda  ^{-1-\frac1{4m}\;(+\varepsilon
)}),
\tag3.22
$$
with $\varepsilon =0$ if $\sigma -\frac14\notin\Bbb Z$. 

2) There exist particular choices of $P_1$ where one has a better
expansion than (3.21):
$$
\Tr(AQ_{1,\lambda ,+})
=\sum_{0\le j <n+\sigma +\frac14}c_{j}\,(-\lambda )
^{\frac{n+\sigma -j}m-1}+(c'_0\log(-\lambda
)+c''_0)(-\lambda )^{-1}+O(\lambda  ^{-1-\frac1{4m}\;(+\varepsilon )}),
\tag3.23
$$
with $\varepsilon =0$ if $\sigma -\frac14\notin\Bbb Z$. 

Then an expansion (3.23) is obtained for general choices of $P_1$ by
use of (3.22).

For point 1) in this program, let us denote
$$
\Cal Q_\lambda =Q_{1,\lambda }-Q_{2,\lambda },\text{ with symbol
}\frak q(x,\xi ,\lambda ).
\tag3.24$$
Then we can write, since $A=P_++G$,
$$
\aligned
A(Q_{1,\lambda} -Q_{2,\lambda })_+&=A\Cal Q_{\lambda ,+}
= (P\Cal Q_\lambda )_++\Cal G_\lambda ,\text{ with }
\\
\Cal G_\lambda &=
-G^+(P)G^{-}(\Cal Q_\lambda )+G\Cal Q_{\lambda ,+}.
\endaligned\tag3.25
$$ 
The last identity refers to a localized situation: In $\Bbb R^n$,
$G^+(P)=r^+Pe^-J$ and $G^-(P)=Jr^-Pe^+$, where $e^\pm$ denote
extension by zero from $\rnpm$ to $\Rn$, $r^\pm$ denote restriction
from $\Rn$ to $\rnpm$, and $J$ maps $u(x',x_n)$ to $u(x',-x_n)$,
cf.\ \cite{G1, p.\ 252 and (A.32)}. In the present situation,
$P_+\Cal Q_{\lambda 
,+}=(P\Cal Q_\lambda )_+-G^+(P)G^-(\Cal Q_\lambda )$. 

The desired formula for the $\psi $do term can be found pointwise in
$x\in \widetilde X$ by use
of Theorem 2.2, and then integrated over $X$. It is the singular
Green term that requires a new
effort.  

\proclaim{Theorem 3.6} Let $A=P_++G$, of order $\sigma \in\Bbb R$ with $G$ of
class 0, assuming $P=0$ if $\sigma \notin\Bbb Z$. Let $P_1$ and $P_2$
be auxiliary elliptic differential operators of order $m$, as
described in the beginning of this section, with $m>n+\sigma $.

The singular Green part $\Cal G_\lambda $  of  $A(Q_{1,\lambda
}-Q_{2,\lambda })$  is
of order $\sigma -m$, class $0$ and regularity $\sigma +\frac12$. 
Consequently, in local coordinates, its normal trace $\Cal S_\lambda
=\tr_n\Cal G_\lambda $
is a $\psi $do on $\Bbb R^{n-1}$ of order $\sigma -m$ and regularity
$\sigma +\frac14$. Denoting its symbol $\frak s(x',\xi ',\lambda
)\sim\sum_{l\ge 0}\frak s_{\sigma -m-l}(x',\xi ',\lambda )$, we have
the trace expansion 
$$
\Tr\Cal G_\lambda 
=\sum_{1\le j < n+\sigma +\frac14} \frak s_{j}\,(-\lambda )
^{\frac{n+\sigma -j}m-1}+O(\lambda  ^{-1-\frac1{4m}\;(+\varepsilon )}),
\tag3.26$$
with $\varepsilon =0$ if $\sigma -\frac14\notin\Bbb Z$, $\varepsilon
>0$ if  $\sigma -\frac14\in\Bbb Z$. Here
$$
\gathered
\frak s_j=\tilde {\frak s}_{j-1}=\int\tr \tilde{\frak
s}_{j-1}(x')\,dx',\text{ where}\\
\tilde {\frak s}_{l}(x')=
\int_{\Bbb R^{n-1}}{\frak s}^h_{\sigma -m-l}(x',\xi ',-1)\,\d\xi'.
\endgathered
\tag3.27
$$ 
 
\endproclaim

\demo{Proof} In a localized situation, $\Cal G_\lambda $ consists of the
operators $-G^+(P)G^-(\Cal Q_\lambda )$ 
and $G\Cal Q_{\lambda ,+}$.
We first study $G^+(P)G^-(\Cal Q_{\lambda }) $. 

As shown in Section 2 (where $Q_{1,\lambda },Q_{2,\lambda }$
were denoted $Q_\lambda ,Q'_\lambda $), the symbol $\frak q$ of $\Cal
Q_\lambda $ in (3.24) has an expansion $\sum_{j\in\Bbb N}\frak q_{-m-j}$ 
in homogeneous symbols $\frak q_{-m-j}$ that are sums
of terms of the form generalizing
(1.5) with at least two factors
$(p_{1,m}-\lambda )^{-1} $ or $(p_{2,m}-\lambda )^{-1} $. 
Since $P_1$
and $P_2$ are differential operators, we need not smooth out around
$\xi =0$, but can take the exact symbols. Since
$p_{1,m}$ and $p_{2,m}$ are scalar, the factors
$q_{i,-m}=(p_{i,m}-\lambda )^{-1}$ can be collected to the
right, so in fact the terms in $\frak q$ are of the form
$$
f(x,\xi ,\lambda )=f_0(x,\xi )\,q_{1,-m}^{\nu _1}\,q_{2,-m}^{\nu _2},\tag3.28
$$ 
with $\nu _1+\nu _2\ge 2$ and $f_0$ polynomial in $\xi $.
Then for each $j\ge 0$, $\frak q_{-m-j}$ is a sum of terms of the form
$r'(x,\xi )q''(x,\xi ,\lambda )$, where $r'$ is the symbol of a
differential operator of order $m$ independent of $\lambda $ and
$q''$ is of order $-2m-j$, 
likewise with structure as in (3.28), smooth in all
variables (for $|\xi |+|\lambda |\ne 0$). The operator $\op(r'q'')$ can
be further decomposed into a finite sum of terms $RQ'=\op(r)\op(q')$, where 
$r$ and $q'$ have a similar structure as $r'$ and $q''$ (we 
need this modification to get a composition of two operators instead of a
product of symbols). Now we treat each term
$$  
G^+(P)G^-(RQ')
$$
separately.
We are working in $\Bbb R^n$ where the manifold corresponds to $\crnp$,
and can assume that the symbols of $P$, $R$ and $Q'$ are defined on
$\Bbb R^n$. Here we write:
$$
\multline
G^-(RQ')=Jr^-RQ'e^+=Jr^-R(e^+r^++e^-r^-)Q'e^+\\
=Jr^-Re^-JJr^-Q'e^+=
\overline R_+G^-(Q'),\endmultline\tag3.29
$$
where we have used that $r^-Re^+=0$ since $R$ is a differential
operator, and denoted $JRJ=\overline R$, again a differential
operator.
Thus $
G^+(P)G^-(RQ')=G^+(P)\overline R_+G^-(Q')$,
where $G^+(P)\overline R_+$ is a $\lambda $-independent s.g.o.\ of
order $\sigma +m$ and class $ m$. It enters in the
parameter-dependent calculus as an operator of order $\sigma +m$,
class $m$ and regularity $\sigma +\frac12$, cf.\ \cite{G1, (2.3.55)}.
(It is the presence of the normal derivatives of order $\le m$ in the
differential operator $\overline R$ that
brings the regularity down to $\sigma +\frac12$, not $\sigma +m$ as
in the considerations for closed manifolds, but the gain of
$\frac12$ will be just enough to serve our purposes.)
Composing with $G^-(Q')$ of order $-2m-j$, class 0 and
regularity $+\infty $, we find that$$
\multline
G^+(P)G^-(RQ')=G^+(P)\overline R_+G^-(Q')\\
\text{ is of order $\sigma
-m-j$, class 0 and regularity }\sigma +\tfrac12. 
\endmultline\tag3.30
$$ 

Collecting the terms (finitely many for each order) we find that the
homogeneous terms in $\frak q$ contribute to an s.g.o.\ of order
$\sigma -m$, class 0 and regularity $\sigma +\frac12$. Since the
remainder of $\frak q$ after subtraction of $N$ homogeneous terms is
$O(\ang{\xi ,\mu }^{-m-N-1})$, its contribution will, when $N$ gets
large, reach arbitrary low orders and estimates $O(|\lambda
|^{-N'})$ for any $N'$, so it complies with the regularity $\sigma
+\frac12$.   

There is a very similar proof for $G\Cal Q_{\lambda ,+}$. Again we
use that each $\op(\frak q_{-m-j})$ can be written as a finite sum of
terms $RQ'$, where $R$ is a differential operator of order $m$ and
$Q'$ has symbol structure as in (3.28) and order $-2m-j$. Now for
each term, since $G^+(R)=0$,
$$
G(RQ')_+=GR_+Q'_+,
$$
where $GR_+$ is a parameter-independent s.g.o.\ of order $\sigma +m-j$
and class $m$, hence has regularity $\sigma +\frac12$ when taken into
the parameter-dependent theory. Then $GR_+Q'_+$ is of order $\sigma
-m-j$, class 0 and regularity $\sigma +\frac12$.
Collecting the terms and treating remainders as above, we get that
$G\Cal
Q_{\lambda ,+}$ has order $\sigma -m$, class 0 and regularity $\sigma
+\frac12$.

This shows the asserted symbol properties of $\Cal G_\lambda $.
Its normal trace $\Cal S_\lambda $ is of order $\sigma -m$ and
regularity $\sigma 
+\frac14$ by Lemma 3.4.
\comment
Let us denote the s.g.o.\ part by $\Cal G_\lambda $, with normal
trace $\Cal S_\lambda $:
$$
\Cal G_\lambda =-G^+(P)G^-(\Cal Q_\lambda )+G\Cal Q_{\lambda ,+},\quad
\Cal S_\lambda =\tr_n \Cal G_\lambda ,\tag 3.37
$$
with symbols $\frak g(x',\xi ',\xi _n,\eta _n,\lambda )$ resp.\ $\frak
s(x',\xi ',\lambda )$. By Theorem 3.6, $\Cal S_\lambda $ is of order $\sigma
-m$ and regularity $\sigma +\frac14$. 
\endcomment
By Lemma
3.5, its kernel 
has an expansion on the diagonal:
$$
\aligned
K(\Cal S_\lambda ,x',x')&=\sum_{0\le l\le n-1+\sigma }\tilde
{\frak s}_{l}(x')(-\lambda )
^{\frac{n-1+\sigma -m-l}m}+O(\lambda ^{\frac{\sigma -m-\sigma
-\frac14\,(+\varepsilon '')}{ m
}})\\
&=\sum_{0\le l\le n-1+\sigma }\tilde
{\frak s}_{l}(x')(-\lambda )
^{\frac{n-1+\sigma -l}m-1}+O(\lambda ^{-1 -\frac1{4 m }\,(+\varepsilon )}),\endaligned\tag3.31
$$
with $\varepsilon =0$ unless $\sigma +\frac14\in\Bbb Z$, and 
$\tilde {\frak s}_{l}(x')$ defined as in (3.27).
In the proof, the lemma is applied for each ray; the ray $
\rmi$ gives the value (3.27) for the coefficients, and the holomorphy
assures that their values are the same on the other rays
(as in \cite{GS1, Lemma 2.3}).
Finally, integration in $x'$ of the fiber trace then gives:
$$
\aligned
\Tr_{\rnp}\Cal G_\lambda &= \Tr_{\Bbb R^{n-1}}\Cal S_\lambda 
=\sum_{0\le l < n-1+\sigma +\frac14} \tilde {\frak s}_{l}\,(-\lambda )
^{\frac{n-1+\sigma -l}m-1}+O(\lambda  ^{-1-\frac1{4m}
\;(+\varepsilon )})\\
&=\sum_{1\le j < n+\sigma +\frac14} \frak s_{j}\,(-\lambda )
^{\frac{n+\sigma -j}m-1}+O(\lambda  ^{-1-\frac1{4m}\;(+\varepsilon )}),
\endaligned\tag3.32
$$
with $\frak s_j$ defined as in (3.27).
\qed
\enddemo 

Observe a direct consequence:

\proclaim{Corollary 3.7} Assumptions as in Theorem {\rm 3.6}.

The trace of $A(Q_{1,\lambda }-Q_{2,\lambda })_+$ has an expansion
{\rm (3.22)}.
\endproclaim 

\demo{Proof} If $\sigma \notin\Bbb Z$, there is no $\psi $do part,
and the expansion is (3.26). If $\sigma \in\Bbb
Z$, the $\psi $do part has an expansion
$$
\Tr\bigl((P(Q_{1,\lambda }-Q_{2,\lambda }))_+\bigr)
=\sum_{0\le j \le n+\sigma }c_{j,+}\,(-\lambda )
^{\frac{n+\sigma -j}m-1}+O(\lambda  ^{-2+\varepsilon ' })
\tag3.33
$$
(any $\varepsilon ' >0$), found from (2.9) by taking fiber traces
and integrating over
$X$. 
When we add this to (3.26), we find (3.22).
\qed
\enddemo 

Now we turn to point 2) in the program for showing (3.23) in general.

\proclaim{Lemma 3.8} Let $P_0$ be selfadjoint positive of order $2 $
with scalar
principal symbol, and let $A=P_++G$ be as above. For $k$ so large that
$2k>n+\sigma $, there is a trace expansion for $\lambda \to\infty $
in $\Bbb C\setminus \rp$:
$$
\Tr(A(P_0^k-\lambda )^{-1}_+)\sim  \sum_{j\ge
0}\!c_{j}(-\lambda)
^{\frac{n+\nu - j}{2k}-1}
+\sum_{l\ge
0}( c'_l\log (-\lambda )+c''_l)(-\lambda )
^{-\frac l{2k}-1}. \tag3.34
$$   
Here $c_0'=\frac 1{2k}\operatorname{res}A$.
\endproclaim 

\demo{Proof} We here use (3.4) for $P_0$, translating it to a
statement on the meromorphic structure of the generalized zeta
function $\zeta (A,P_{0,+},s)$,
which allows replacing $P_0$ by $P_0^k$; this gives the structure of
$\zeta (A, (P_0^k)_+,s)$, which translates back to a trace expansion
(3.34). Here \cite{GS2, Prop.\ 2.9, Cor.\ 2.10 and Cor.\ 3.5} are
used. In details:

We define $\zeta (A,P_{0,+},s)$ and $\zeta (A,(P_0^k)_+,s)$ as the
meromorphic extensions of \linebreak $\Tr(A(P_0^{-s})_+)$ resp.\
$\Tr(A(P_0^{-sk})_+)$, defined \`a priori for large
$\operatorname{Re}s$.  It is well-known that the expansion (3.4)
implies the following meromorphic structure of $\zeta (A, P_{0,+},
s)$: 
$$
\Gamma (s)\zeta (A, P_{0,+}, s)\sim \sum_{j\ge 0}\frac{\tilde
c_j } {s+\frac{j-n-\sigma }{2}} +\sum_{l\ge 0}\Bigl(\frac{\tilde
c'_l}{(s+\frac l2)^2} +\frac{\tilde c''_l} {s+\frac l2}\Bigr)
\tag3.35$$ 
(by
use of e.g.\ \cite{GS2, Cor.\ 2.10}).  Dividing out the Gamma factor,
we obtain a meromorphic structure somewhat similar to (3.35), 
$$
\zeta (A, P_{0,+}, s)\sim \sum_{j\ge 0}\frac{\tilde
a_j } {s+\frac{j-n-\sigma }{2}} +\sum_{l\ge 0}\Bigl(\frac{\tilde
a'_l}{(s+\frac l2)^2} +\frac{\tilde a''_l} {s+\frac l2}\Bigr)
\tag 3.36
$$ 
{\it except that the double poles vanish for $l$ even}, since they are
turned into simple poles by the cancellations from the zeros of
$\Gamma (s)^{-1}$ at 0, $-1$, $-2$, \dots. Since $P_0$ is selfadjoint
positive, the complex powers agree with the definition by spectral
theory, so $P_0^{-s}=(P_0^k)^{-s'}$, $s'=s/k$. Then we can replace
the formula for $P_0^{-s}$ by the formula for 
$(P_0^k)^{-s'}$ simply by  replacing the variable $s$ by $s'k$, so we get
$$
\aligned
\zeta (A, (P_0^k)_+, s')&\sim \sum_{j\ge 0}\frac{\tilde
a_j } {s'k+\frac{j-n-\sigma }{2}} +\sum_{l\ge 0}\Bigl(\frac{\tilde
a'_l}{(s'k+\frac l2)^2} +\frac{\tilde a''_l} {s'k+\frac l2}\Bigr)\\
&\sim \sum_{j\ge 0}\frac{\tilde
b_j } {s'+\frac{j-n-\sigma }{2k}} +\sum_{l\ge 0}\Bigl(\frac{\tilde
b'_l}{(s'+\frac l{2k})^2} +\frac{\tilde b''_l} {s'+\frac l{2k}}\Bigr),
\endaligned
$$ 
with the double poles vanishing for $l$ even. 
Multiplication by
$\Gamma (s')$ gives still another expansion
$$
\Gamma (s')\zeta (A, (P_0^k)_+, s')\sim \sum_{j\ge 0}\frac{\tilde
d_j } {s'+\frac{j-n-\sigma }{2k}} +\sum_{l\ge 0}\Bigl(\frac{\tilde
d'_l}{(s'+\frac l{2k})^2} +\frac{\tilde d''_l} {s'+\frac l{2k}}\Bigr)
\tag 3.37$$
where we get double poles back at the values where $l/k$ is even (a
subset of the set where they were removed before).

Finally, we use \cite{GS2, Prop.\ 2.9} in the direction from $\zeta
(s)$ to $f(\lambda )$, in 
the same way as in the proof of \cite{GS2, Cor.\ 3.5}.
The cited proposition shows how the meromorphic structure of $\Gamma
(1-s')\Gamma (s')\zeta (A, (P_0^k)_+, s')$ carries over to
an asymptotic expansion of $f(-\lambda
)=$ \linebreak$\Tr(A(P_0^k-\lambda )_+^{-1})$. The needed exponential decrease for
$|\operatorname{Im}s'|\to \infty $ follows from the similar property
of $\Gamma
(1-s)\Gamma (s)\zeta (A, P_{0,+}, s)$. That
$f(-\lambda )$  satisfies an $O(|\lambda |^{-\alpha })$ estimate
(with $\alpha >0$) for
$\lambda \to \infty $ in the considered sector is assured by
(3.21) above, with $m=2k>n+\sigma $. The positivity of $P_0$ assures
that $f$ is
regular at 0. The method introduces
some possible new integer poles on the positive real axis (coming from
$\Gamma (1-s')$), but in the end result they are not present, since we
already have the corresponding part of the expansion in powers known
from (3.21). 

It is known from \cite{GSc1} that $\tilde
c'_0=\frac12\operatorname{res}A$ in (3.35), also equal to the coefficient
of $(\lambda )^{-1}\log(-\lambda )$ in the corresponding resolvent
trace expansion. Following the reduction, we see that $\tilde d'_0$
in (3.37) equals $\frac1k \tilde c'_0=\frac1{2k}\operatorname{res}A$. 
\qed 
\enddemo 

\proclaim{Theorem 3.9} Assumptions as in Theorem {\rm 3.6}. 
 
$\Tr(AQ_{1,\lambda ,+})$ has
a trace expansion {\rm (3.23)}; in particular, $C_0(A,P_{1,+})$ is
well-defined as the coefficient of $(-\lambda )^{-1}$, and
$c'_0=\frac1m \operatorname{res}A$.
\endproclaim 

\demo{Proof}
First let $m$ be even $=2k$. Then we can compare an arbitrary
auxiliary operator $P_1$ with $P_3=P_0^k$ from Lemma 3.8 (with resolvent $Q_{3,\lambda
}=(P_3-\lambda )^{-1}$). Here (3.22) for the trace difference and
(3.34) for $\Tr(AQ_{3,\lambda ,+})$ add up to give:
$$
\multline
\Tr(AQ_{1,\lambda ,+})=\Tr\bigl(A(Q_{1,\lambda }-Q_{3,\lambda
})_+\bigr)+\Tr(AQ_{3,\lambda ,+})\\
=\sum_{0\le j < n+\sigma +\frac14 }c_{j}\,(-\lambda )
^{\frac{n+\sigma -j}m-1}+(c'_0\log(-\lambda
)+c''_0)(-\lambda )^{-1}+O(\lambda  ^{-1-\frac1{4m}\;(+\varepsilon )}),
\endmultline\tag3.38
$$
with $c'_0=\frac1m\operatorname{res}A$.
So the assertion holds for $m$ even.

Next, let $m$ be odd. Necessarily, $P_1$ cannot have its spectrum in a
sector with opening $<\pi $, since the principal symbol is odd in
$\xi $, so iterated powers are not easy to use (e.g., for a
selfadjoint Dirac
operator $D$, $(D^2)^\frac12=|D|$ is different from $D$). Instead we
shall use an idea of doubling up, found in \cite{GS1}. For a given
$P_1$ of order $m$, consider$$
\Cal P_1=\pmatrix 0 & -P_1^*\\ P_1&0 \endpmatrix,\quad \Cal A=\pmatrix
A & 0\\ 0&A  \endpmatrix,
$$
acting in the bundle $E\oplus E$. $\Cal P_1$ is skew-selfadjoint, with
resolvent $$
(\Cal P_1-\lambda )^{-1}=\pmatrix -\lambda (P_1^*P_1+\lambda ^2)^{-1}
& P_1^*(P_1P_1^*+\lambda ^2)^{-1}\\ -P_1(P_1^*P_1+\lambda
^2)^{-1}&-\lambda (P_1P_1^*+\lambda ^2)^{-1} \endpmatrix,\quad 
$$
for $\lambda \in\Bbb C\setminus i\Bbb R$. Now$$
\Tr(\Cal A(\Cal P_1-\lambda )_+^{-1})=-\lambda \Tr(A(P_1^*P_1+\lambda
^2)_+^{-1})-\lambda \Tr(A(P_1P_1^*+\lambda ^2)_+^{-1}).
$$
In the right-hand side, $P_1^*P_1$ and $P_1P_1^*$ are sefadjoint
elliptic of even order $2m$ (and $\ge 0$), so by the result already
shown for even-order auxiliary operators, applied to the two traces,
we get:
$$
\multline
\Tr(\Cal A(\Cal P_1-\lambda )_+^{-1})=-\lambda \Bigl(\sum_{0\le j < n+\sigma +\frac14 }a_{j}\,\lambda ^{2(
\frac{n+\sigma -j}{2m}-1)}+(a'_0\log(\lambda
^2)+a''_0)\lambda ^{-2}\\
\quad+O(\lambda  ^{2(-1-\frac1{8m}\;(+\varepsilon
)}))
\Bigr)\\
= \sum_{0\le j < n+\sigma +\frac14 }b_{j}\,(-\lambda )
^{\frac{n+\sigma -j}{m}-1}+(2a'_0\log(-\lambda
)+a''_0)(-\lambda )^{-1}+O(\lambda  ^{-1-\frac1{4m}\;(+\varepsilon 
)})
,\endmultline\tag3.39
$$
with coefficients modified because of powers of $-1=e^{i\pi }$; here
$a'_0=2\frac1{2m}\operatorname{res}A=\frac1m\operatorname{res}A$.
Now $\Cal P_1$ can be compared with $$
\Cal P_2=\pmatrix P_1 & 0\\ 0&P_1 \endpmatrix,
$$
and a calculation as in (3.38) gives that 
$\Tr(\Cal A(\Cal P_2-\lambda )_+^{-1})$ likewise has an expansion as in
the last line of (3.39), with the same coefficient $2a'_0$ of the
logarithmic term. Then 
$\Tr (A(P_1-\lambda )^{-1}_+)=\frac12 \Tr(\Cal A(\Cal P_2-\lambda
)_+^{-1})$ likewise has an expansion, with log-coefficient $ a'_0$.\qed
\enddemo 

One can also see from these proofs that the value of $C_0(A,P_{1,+})$
{\it modulo 
local terms} is as described in \cite{GSc2}, namely, in local
coordinates, a sum of integrals over $X$ resp.\ $X'$ of finite part
integrals in $\xi $ resp.\ $\xi '$ of the symbols of $P$ resp.\
$\tr_n G$.

Now the coefficient of $(-\lambda )^{-1}$  in $\Tr(A(Q_{1,\lambda }-
Q_{2,\lambda })_+)$ will be studied in detail.

\comment
Observe a direct consequence:

\proclaim{Corollary 3.8} Assumptions as in Theorem {\rm 3.6}.

The trace of $A(Q_{1,\lambda }-Q_{2,\lambda })_+$ has an expansion 
$$
\Tr\bigl(A(Q_{1,\lambda }-Q_{2,\lambda })_+\bigr)
=\sum_{0\le j < n+\sigma +\frac14}d_{j}\,(-\lambda )
^{\frac{n+\sigma -j}m-1}+O(\lambda  ^{-1-\frac1{4m}\; (+\varepsilon )}),
\tag3.37
$$
where $\varepsilon =0$ unless $\sigma +\frac14$ is integer.

Moreover, when $m$ is even, $\Tr AQ_{1,\lambda ,+}$ has an
expansion as in {\rm (3.23)}, so that $C_0(A,P_{1,+})$ is well-defined.
\endproclaim 

\demo{Proof} For the $\psi $do part (entering only if $\sigma \in\Bbb
Z$) we have an expansion
$$
\Tr\bigl((P(Q_{1,\lambda }-Q_{2,\lambda }))_+\bigr)
=\sum_{0\le j < n+\sigma +\frac14 }c_{j,+}\,(-\lambda )
^{\frac{n+\sigma -j}m-1}+O(\lambda  ^{-2+\varepsilon ' })
\tag3.33
$$
(any $\varepsilon ' >0$), found from (2.9) by taking fiber traces
and integrating over
$X$. 
When we add this to (3.26), we find (3.37).

For $m$ even, we now have all the ingredients to conclude as in (3.23).
\qed
\enddemo

It is seen from Corollary 3.7 that regardless of whether the
operators $AQ_{i,\lambda ,+}$ have trace expansions down to lower
powers than
 $(-\lambda )^{-1}$ (we believe that they have so, but
have only verified it for even $m$), the difference between them
certainly has such a trace expansion. So we can {\it
define} the trace defect in any case by$$
\operatorname{def}(A,P_{1,+}, P_{2,+})=\text{ coefficient of
}(-\lambda )^{-1} \text{ in (3.37);}
\tag3.39
$$
it equals $C_0(A,P_{1,+})-C_0(A,P_{2,+})$ when both terms are
defined (in particular when $m$ is
even).
\endcomment

Note that when $n+\sigma \in\Bbb N$, the sum in (3.22)
goes from $0$ 
to $n+\sigma $ and the last term is $d_{n+\sigma }(-\lambda )^{-1}$.
When $n+\sigma \notin \Bbb N$, we see that there is no term with
$(-\lambda )^{-1}$ in the expansion, so$$
C_0(A,P_{1,+})-C_0(A,P_{2,+})=0\text{ if }n+\sigma
\notin\Bbb N.\tag3.40 
$$
We shall finally show:

\proclaim{Theorem 3.10} Assumptions as in Theorem {\rm 3.6}. One has that
$$
C_0(A,P_{1,+})-C_0(A,P_{2,+})=-\tfrac 1m \operatorname{res}(A(\log
P_1-\log P_2)_+).\tag 3.41
$$
\endproclaim 

\demo{Proof} 
Denote
$$
 L=\log P_1-\log P_2,\tag3.42
$$
with symbol $l(x,\xi )$; in view of (1.14), it is classical of order 0, 
and (cf.\ (3.24))
the homogeneous terms $l_{-j}(x,\xi )$ are
determined for $|\xi |\ge 1$ by the formulas$$
l_{-j}(x,\xi )=\tfrac i{2\pi }\int_{\Cal C'}\log\lambda
\,\frak q_{-m-j}(x,\xi ,\lambda )\, d\lambda ,\tag3.43
$$
where $\Cal C'$ is a closed curve in $\Bbb C\setminus\crm$ encircling
the values of 
$p_{1,m}(x,\xi )$ and $p_{2,m}(x,\xi )$. From the fact that $P_1$ and $P_2$ are differential
operators, it is easily checked that $L$ satisfies the transmission
condition at $x_n=0$.

We have that $$  
A(\log P_1-\log
P_2)_+=AL_+=(P_++G)L_+=(PL)_+-G^+(P)G^-(L)+GL_+.  
\tag3.44$$ 
According to Fedosov, Golse, Leichtnam and Schrohe
\cite{FGLS} (with the sign of the s.g.o.-term corrected in
\cite{GSc1}), the residue is determined by the formula
$$\multline
\operatorname{res}(AL_+)=\int_{\rnp}\int_{|\xi
|=1}\tr\operatorname{symb}_{-n} (PL)\,\d S(\xi )dx\\ +\int_{\Bbb
R^{n-1}}\int_{|\xi' |=1}\tr\operatorname{symb}_{1-n}
(\tr_n(-G^+(P)G^-(L)+GL_+))\,\d S(\xi ')dx', 
\endmultline\tag3.45 $$
where $\operatorname{symb}_k$ stands for ``the homogeneous term of
degree $k$ in the symbol of the operator''.

Consider first the case where $\sigma \notin\Bbb Z$, $P=0$. 
Then the left-hand side in (3.41) is zero in view of (3.40), and the
right-hand side is zero, since $A(\log P_1-\log P_2)_+$ is an s.g.o.\
of noninteger order. So the formula
is verified for $\sigma \notin\Bbb Z$, and we can restrict the
attention to the
case where $\sigma \in\Bbb Z$.

The calculations leading
to Theorem 2.2 show that 
$$ -\tfrac1m\int_{\rnp}\int_{|\xi
|=1}\tr\operatorname{symb}_{-n} (PL)\,\d S(\xi
)dx=\int_{\rnp}\int_{\Rn }\tr
\operatorname{symb}^h_{-m-n} (P\Cal Q_\lambda )|_{\lambda =-1}\,\d
\xi dx; \tag3.46 
$$ this gives the $\psi $do part of the desired
formula.

Now consider the s.g.o.\ part.
For the operator $\Cal G_\lambda $ and its normal trace $\Cal
S_\lambda $, we denote the symbols $\frak g(x',\xi ',\xi _n,\eta _n,\lambda )$ resp.\ $\frak
s(x',\xi ',\lambda )$.
Moreover, we denote$$
G'=-G^+(P)G^-(L)+GL_+,\quad
S'=\tr_n(-G^+(P)G^-(L)+GL_+),
\tag 3.47
$$
with symbols  $g'(x',\xi ',\xi _n,\eta _n)$, $s'(x',\xi ')$.

From (3.27) we have in particular:
$$
\tilde {\frak s}_{\sigma +n-1}(x)=
\int_{\Bbb R^{n-1}}\frak s^h_{ -m+1-n}(x',\xi ',-1)\,\d\xi',
\tag3.48
$$ 
and the integral of its fiber trace gives the contribution to 
$C_0(A,P_{1,+})-C_0(A,P_{2,+})$.

In the following, consider first the case where $\frak q$ is
independent of $x_n$.

 The term of order $1-n$ in the symbol of $S'$ is constructed for
$|\xi '|\ge 1$ as the term of homogeneity degree $1-n$ in the
symbol$$
\aligned
s'(x',\xi ')&=\int_{\Bbb R}(-g^+(p)\circ g^-(l)+g\circ l_+)(x',\xi
',\xi _n,\xi _n)\,\d\xi _n\\
&=\int_{\Bbb R}\bigl(-g^+(p)\circ g^-(\tfrac i{2\pi }\int_{\Cal C''}\log
\lambda \, \frak q\,
d\lambda )+g\circ (\tfrac i{2\pi }\int_{\Cal C''}\log \lambda \,\frak
q\,d\lambda )_+\bigr)\,\d\xi _n;\\
\endaligned
\tag3.49
$$
here $\Cal C''$ is a curve in $\Bbb C\setminus \crm$ formed as the
boundary of a set $V_{r,\theta }$ (1.18) with $\theta $, $r$ and
$\varepsilon $ taken
so small that 
the principal symbols $p_{1,m}(x,\xi )$  and $p_{2,m}(x,\xi )$ take values in
the complement of $ V_{r+\varepsilon ,\theta +\varepsilon }$ for all
$x$, all $|\xi 
'|\ge 1$. Such sets exist since $p_{1,m}$ and
$p_{2,m}$ are homogeneous of degree $m$ in $\xi $ for $|\xi |\ge 1$ and the
ellipticity condition holds uniformly in $x$ (originally running in
the compact manifold $X$). We have:
$$
\aligned
&s'_{1-n}(x',\xi ')\\
&=\sum_{j+k+|\alpha |= \sigma +n-1}\tfrac
{(-i)^{|\alpha |}}{\alpha !}\Bigl(
\int_{\Bbb R}(-\partial _{\xi '}^\alpha g^+(p)_{\sigma -j}\circ _ng^-(\tfrac i{2\pi }\int_{\Cal C''}\log
\lambda \, \partial _{x'}^\alpha \frak q_{-m-k}\,
d\lambda )\d\xi _n\\
&\quad +\int_{\Bbb R} \partial _{\xi '}^\alpha g_{\sigma -j}\circ _n(\tfrac i{2\pi }\int_{\Cal C''}\log
\lambda \,\partial _{x'}^\alpha \frak q_{-m-k}\,d\lambda )_+)\,\d\xi _n\Bigr),
\endaligned
\tag3.50
$$
where $\circ_n$ stands for symbol composition with respect to the
normal variables (cf.\ \cite{G1, Sect.\ 2.6}), and we denote the
homogeneous term
of order $r$ in an s.g.o.\ symbol $g$ by $g_r$ (it is of {\it degree} $r-1$;
this index was used in \cite{G1}).
There are finitely many terms. In each term, the integration in
$\lambda $ and the factor $\log \lambda $ can be moved outside
$\circ_n$ and $g^-$, since these 
operations preserve the holomorphy in $V_{r+\varepsilon ,\theta
+\varepsilon }$ and preserve sufficient decrease in $\lambda $ for $|\lambda
|\to \infty $ (in view of the detailed rules in \cite{G1, Sect. 2.6}
and the analysis in Theorem 3.6). Furthermore, the integrations in $\lambda
$ and $\xi _n$ can be interchanged. So if we define$$
\aligned
\varphi (x',\xi ',\lambda )&=\sum_{j+k+|\alpha |= \sigma +n-1}\tfrac
{(-i)^{|\alpha |}}{\alpha !}
\int_{\Bbb R}\Bigl[-\partial _{\xi '}^\alpha g^+(
p)_{\sigma -j}\circ _ng^-( \partial _{x'}^\alpha \frak q_{-m-k}(x',0,\xi
,\lambda ))\\
&\quad + \partial _{\xi '}^\alpha g_{\sigma -j}\circ _n( \partial _{x'}^\alpha
\frak q_{-m-k}(x',0,\xi ,\lambda  ))_+\Bigr]\,\d\xi _n,
\endaligned
\tag3.51
$$
we have that$$
s'_{1-n}(x',\xi ')=\tfrac i{2\pi }\int_{\Cal C''}\log \lambda \,
\varphi (x',\xi ',\lambda )\,d\lambda.\tag3.52 
$$
An application of Lemma 1.2 gives:
$$
\tfrac i{2\pi }\int_{\Cal C''}\log \lambda \,
\varphi (x',\xi ',\lambda )\,d\lambda =-\int_{-\infty }^0\varphi (x',\xi ',t)\,dt.\tag3.53
$$
One checks from (3.51) that $\varphi $ has the quasi-homogeneity property
$\varphi (x',t\xi ',t^m\lambda )$\linebreak$=t^{-m-n+1}\varphi (x',\xi ',\lambda )$
for $|\xi '|\ge 1$, $t\ge 1$.
 Taking strictly homogeneous
symbols everywhere gives $\varphi ^h$, which is integrable at
$\xi '=0 $ for $\lambda \ne 0$ in view of the regularity properties 
shown in Theorem 3.6.
Now we
can apply Lemma 1.3 with dimension $n$
replaced by $n-1$, finding that
$$
\aligned
&\int_{|\xi '|=1}s'_{1-n}(x',\xi ')\,\d S(\xi ')=-\int_{|\xi
'|=1}\int_{-\infty }^0\varphi ^h(x',r\xi ',t)\,dt\d S(\xi ')\\
&=-m\int_{\Bbb R^{n-1}}\varphi ^h(x',\xi ',-1)\,\d\xi '\\
&=-m\sum_{j+k+|\alpha |= \sigma +n-1}\tfrac
{(-i)^{|\alpha |}}{\alpha !}
\int_{\Bbb R^{n-1} }\int_{\Bbb R}\Bigl[-\partial _{\xi '}^\alpha g^+(
p)^h_{\sigma -j}\circ _ng^-( \partial _{x'}^\alpha \frak
q^h_{-m-k}(x',0,\xi ,-1)\\
&\quad + \partial _{\xi '}^\alpha g^h_{\sigma -j}\circ _n( \partial
_{x'}^\alpha 
\frak q^h_{-m-k}(x',0,\xi ,-1 ))_+\Bigr]\,\d\xi _n\d\xi ',
\endaligned
\tag3.54
$$
which we recognize as
$$
=-m\int_{\Bbb R^{n-1}}\frak s^h_{-m+1-n}(x',\xi ',-1)\,\d\xi '.\tag3.55
$$

This shows that the contribution from $S'=\tr_n G'$ (cf.\ (3.47))
matches the coefficient $\tilde {\frak s}_{n+\sigma }(x')$ of $(-\lambda
)^{-1}$ in the diagonal kernel 
expansion of $\Cal S_\lambda =\tr_n \Cal G_\lambda $, pointwise in
$x'$, cf.\ (3.48), (3.27).
Integration of the fiber trace in $x'$ gives $$
\operatorname{res}(G')=-m\,\tilde{\frak s}_{n-1+\sigma },\tag3.56
$$
where $\tilde{\frak s}_{n-1+\sigma }$ is the coefficient of $(-\lambda
)^{-1}$ in the trace expansion of $\Cal S_\lambda $ (and 
the trace expansion of $\Cal G_\lambda $), cf.\ (3.32).
Adding this identity to (3.46), we find (3.41).

There remains to include the case where the symbol $\frak q$
depends on $x_n$, but this is easy to do. One takes a Taylor
expansion of $\frak q$ in $x_n$ at $x_n=0$; since a factor $x_n^k$ lowers the
order in the resulting s.g.o.s by $k$ steps (cf.\ \cite{G1, Lemma 2.4.3}), only the first
$\sigma +n$ terms can contribute to the constants we are studying.
Each of these terms enters by the standard composition rules in a
very similar way as above, only now one also has to keep track of the
effect of powers $x_n^k$. Again this leads to (3.41).
\qed
\enddemo 

\example{Remark 3.11} The proof shows that the identity (3.41) holds
in a partly localized way, namely, the pseudodifferential
contributions from each side match pointwise in $x\in X$ (before
integration in $x$), and for the singular Green contributions, the
$\psi $do's on $X'$ obtained after taking $\tr_n$ match pointwise in
$x'\in X'$ (before integration in $x'$).
\endexample

\example{Remark 3.12} The identity (3.41) holds also when the $P_i$ are
taken of order $m=2$ as in \cite{GSc2}, which necessitates a replacement of
$Q_{i,\lambda }$ by $Q^N_{i,\lambda }$ for a large enough $N$. For,
writing $Q^N_{i,\lambda }=\frac {\partial _\lambda ^{N-1}}{(N-1)!}Q_{i,\lambda
}$, we see
that the term with $(-\lambda )^{-N}$ in $\Tr(A(Q^N_{1,\lambda
}-Q^N_{2,\lambda })_+)$ is found by integration of
compositions where the symbol terms $\frak q^h_{-m-j}$
are replaced by  $\frac {\partial _\lambda
^{N-1}}{(N-1)!}\frak q^h_{-m-j}$. We just give the argument for the
s.g.o.\ part. The analysis in \cite{GSc2} shows
the needed fall-off in $\lambda $ and integrability in $\xi '$ in
this case.  Since
(with notation as in the proof of Theorem 3.10) 
$$
(-\lambda )^{-1}\int_{\Bbb R^{n-1}}\frak s^h_{-m+1-n}(x',\xi
',-1)\,\d\xi '=\int _{\Bbb R^{n-1}}\frak s^h_{-m+1-n}(x',\xi
',\lambda )\,\d\xi '$$
for $\lambda \in\rmi$, an application of 
$\frac {\partial _\lambda
^{N-1}}{(N-1)!}$ gives for the corresponding function $\frak
s^{(N)h}_{-m+1-n}(x',\xi ',\lambda )$ resulting from insertion of the
$\frac {\partial _\lambda
^{N-1}}{(N-1)!}\frak q^h_{-m-j}$:
$$\aligned
\int _{\Bbb R^{n-1}}\frak s^{(N)h}_{-m+1-n}(x',\xi
',\lambda )\,\d\xi '&=\tfrac {\partial _\lambda
^{N-1}}{(N-1)!}\int _{\Bbb R^{n-1}}\frak s^{h}_{-m+1-n}(x',\xi
',\lambda )\,\d\xi '\\
&=\tfrac {\partial _\lambda
^{N-1}}{(N-1)!}\bigl[(-\lambda )^{-1}\int_{\Bbb R^{n-1}}\frak s^h_{-m+1-n}(x',\xi
',-1)\,\d\xi '\bigr]\\
&=(-\lambda )^{-N}\int_{\Bbb R^{n-1}}\frak s^h_{-m+1-n}(x',\xi
',-1)\,\d\xi ',\endaligned\tag3.57
$$ 
showing that the coefficient of $(-\lambda )^{-N}$ in the expansion
of $\Tr(A(Q^N_{1,\lambda
}-Q^N_{2,\lambda })_+)$ obeys the same
formulas as the coefficient of $(-\lambda )^{-1}$ in Theorem 3.10.
Then it is set in relation to the residue in 
exactly the same way as we did there.
\endexample

\example{Remark 3.13} The assumption on scalarity of the principal
symbols of $P_1$ and $P_2$ was convenient in the proof of Theorem
3.6, but can probably be removed; instead one can study the entries
in the matrix symbol of $\Cal Q_\lambda $ individually, collecting
$\lambda $-independent polynomial factors to the left. One can also
use the methods from the following section.
\endexample

\subhead 4. The second trace defect formula for manifolds with
boundary\endsubhead 

Now some words on possible extensions of the other trace defect formula
(2.6) to the situation of $\psi $dbo's. Here we assume $m>\sigma
+\sigma '+n$
in order to have a trace-class operator $[A,A']Q_{\lambda ,+}$. 
Clearly, 
$$
\Tr
([A,A']Q_{\lambda ,+})=\Tr(A[A',Q _{\lambda ,+}]),\tag 4.1$$ 
so one might strive to show that $-m\,C_0([A,A'],P_{1,+})$
should equal
$$\operatorname{res}(A[A',(\log P_1)_+]).\tag 4.2
$$
But there are several problems with such a formula. The $\psi $do part
of $A[A',(\log P_1)_+]$ is $P[P',\log P_1]$, hence classical in view
of (1.14). But there will in addition be s.g.o.-like elements that are not
covered by existing theories. One is $G^+(\log P_1)=r^+\log P_1
e^-J$, which is {\it not} a standard s.g.o., for example,
$G^+(\log(-\Delta ))$ on $\rnp$ has symbol-kernel \linebreak$c(x_n+y_n)^{-1}e^{-|\xi
'|(x_n+y_n)}$ (for $|\xi '|\ge 1$) with a singularity at
$x_n=y_n=0$. Furthermore, compositions of $(\log P_1)_+$ with $\psi $dbo's
will also contain non-standard terms.

We shall proceed in a different way. Namely, we show
for the singular Green part
$\Cal G_\lambda $ of $A[A',Q 
_{\lambda +}]$ that its normal trace $\Cal
S_\lambda $ has
sufficiently good symbol estimates to allow a ``log-transform''
(integration together with $\log \lambda $ over a curve $\Cal C''$ as
in Theorem 3.10) resulting in a classical $\psi $do $S$ over $X'$, such
that the contribution from $\Cal G_\lambda $ equals $-\frac 1m
\operatorname{res}S$. 

As in \cite{GSc2}, we assume that
the $\psi $do's $P$ and $P'$ are of normal order 0.
(Normal order $k$ means that 
the symbol and its derivatives are  
 $O(\ang{\xi _n}^k)$ at the boundary, here $k\le$ the order.
In general, when $P$ satisfies the
transmission condition, it is the sum of 
a $\psi $do of normal order $-1$, a differential operator,
and a $\psi $do vanishing to a very high order at the boundary.)

There is
a delicate argument in \cite{GSc1,2} for showing that terms containing
compositions with $G^{\pm}(Q_\lambda )$ contribute to $C_0$ with
local coefficients; this relies on the exact structure of the symbol
of $Q_\lambda $ at $x_n=0$ as a function of the roots of the
polynomial $p_{1,2}(x',0,\xi ',\xi
_n)-\lambda $ in $\xi _n$. We shall here replace this argument 
with
an argument using that $$
Q_\lambda =(P_1-\lambda )^{-1}=-\tfrac 1{\lambda }+\tfrac1\lambda 
P_1(P_1-\lambda )^{-1}=-\tfrac 1{\lambda }+\tfrac1\lambda 
P_1Q_\lambda ,\tag4.3
$$
where the contributions from the term $\frac1\lambda $ cancel
out in the calculations of commutators. 
A difficulty in using this is that $P_1Q_\lambda $ is only of order
0, not of large negative order. 

It may be remarked that the
calculations in the following do not need scalarity of the principal
symbol of $P_1$, since they do not appeal to commutation of factors
in $q$.

We work in a localized situation (as in Section 3).
The singular Green terms appearing in the treatment of $A[A',Q 
_{\lambda +}]$ are calculated in the
following lemma:

\proclaim{Lemma 4.1} Let $A=P_++G$, $A'=P'_++G'$ of orders $\sigma $ resp.\
$\sigma '$, the $\psi $do's being of normal order $0$  and the s.g.o.s
being of class $0$; assume that $P$ and $P'$ are zero if $\sigma $ or
$\sigma '$ is noninteger. The singular Green part $\Cal G_\lambda $
of $A[A',Q _{\lambda +}]$ is the sum of terms
$$
\Cal G_\lambda =G[G',Q_{\lambda ,+}]+P_+[G',Q_{\lambda
,+}]+G[P'_+,Q_{\lambda ,+}]+G_{1,\lambda },\tag 4.4
$$
 with the following properties:

$1^\circ$ $G[G',Q_{\lambda ,+}]$, $P_+[G',Q_{\lambda ,+}]$ and
$G[P'_+,Q_{\lambda ,+}]$ are singular Green operators
satisfying the primary formulas$$
\aligned
G[G',Q_{\lambda ,+}]&=GG'Q_{\lambda
,+}-GQ_{\lambda ,+}G' ,\\
P_+[G',Q_{\lambda ,+}]&=P_+G'Q_{\lambda
,+}-P_+Q_{\lambda ,+}G',\\
G[P'_+,Q_{\lambda ,+}]&=GP'_+Q_{\lambda
,+}-GQ_{\lambda ,+}P'_+ .
\endaligned\tag4.5
$$
and the secondary formulas
$$
\aligned
G[G',Q_{\lambda ,+}]&=\tfrac 1\lambda GG'P_{1,+}Q_{\lambda
,+}-\tfrac1\lambda GP_{1,+}Q_{\lambda ,+}G',\\
P_+[G',Q_{\lambda ,+}]&=\tfrac 1\lambda P_+G'P_{1,+}Q_{\lambda
,+}-\tfrac1\lambda P_+P_{1,+}Q_{\lambda ,+}G',\\
G[P'_+,Q_{\lambda ,+}]&=\tfrac 1\lambda GP'_+P_{1,+}Q_{\lambda
,+}-\tfrac1\lambda GP_{1,+}Q_{\lambda ,+}P'_+.
\endaligned\tag4.6
$$

$2^\circ$ $P_+[P'_+,Q_{\lambda ,+}]$ is the sum of the $\psi $do term
$(P[P',Q_{\lambda }])_+$ and a singular Green term $G_{1,\lambda }$ satisfying
primarily$$
G_{1,\lambda }=
-G^+(P)G^-([P',Q_\lambda ])
-P_+G^+(P')G^-(Q_\lambda )+P_+G^+(Q_\lambda )G^-(P'),\tag4.7 
$$
with
$$\aligned
G^-([P',Q_\lambda ])&=G^-(P'Q_\lambda )-G^-(Q_\lambda P')\\
&=G^-(P')Q_{\lambda ,+}+\overline{P'}_+G^-(Q_\lambda )
 -G^-(Q_\lambda )P'_+-(\overline{Q_\lambda })_+G^-(P'), 
\endaligned\tag4.8 $$
secondarily
$$
\multline
G_{1,\lambda }
=-\tfrac1\lambda G^+(P)G^-([P',P_1Q_\lambda
])\\
-\tfrac1\lambda P_+G^+(P')\overline P_{1,+}G^-(Q_\lambda
)+\tfrac1\lambda P_+P_{1,+}G^+(Q_\lambda
)G^-(P'),
\endmultline\tag4.9
$$
with$$\aligned
G^-([P',P_1Q_\lambda ])&=G^-(P'P_1Q_\lambda )-G^-(P_1Q_\lambda P')\\
&=G^-(P'P_1)Q_{\lambda ,+}+(\overline {PP_1})_+G^-(Q_\lambda )\\
&\quad-P_{1,+}G^-(Q_\lambda
)P'_{+}-
P_{1,+}(\overline{Q_\lambda })_+G^-(P' ).
\endaligned\tag4.10$$
\endproclaim 

\demo{Proof} The cases in $1^\circ$ follow easily by insertion of (4.3),
since multiplication by $\frac1\lambda$ commutes with $G'$ and with
$P'_+$, and (since $P_1$ is a differential operator)$$
(P_1Q_\lambda )_+=r^+P_1Q_{\lambda }e^+=r^+P_1e^+r^+Q_{\lambda
}e^+=P_{1,+}Q_{\lambda ,+}.\tag4.11
$$

For case $2^\circ$, we calculate as 
follows:$$
P_+[P'_+,Q_{\lambda ,+}]=P_+[P',Q_\lambda
]_+-P_+G^+(P')G^-(Q_\lambda )+P_+G^+(Q_\lambda )G^-(P'),\tag4.12
$$
where we have used in the last expression that $P'$ has normal order
$\le 0$. 
Here, since $G^\pm(\frac1\lambda )=0$,
$$
P_+G^+(P')G^-(Q_\lambda )=\tfrac1\lambda P_+G^+(P')G^-(P_1Q_\lambda )
=\tfrac1\lambda P_+G^+(P')\overline P_{1,+}G^-(Q_\lambda )
$$
as in (3.29); similarly, $$
P_+G^+(Q_\lambda )G^-(P')=\tfrac1\lambda P_+G^+(P_1Q_\lambda
)G^-(P')=
\tfrac1\lambda P_+P_{1,+}G^+(Q_\lambda )G^-(P')
$$
in view of (4.11). 
This explains the last two terms in (4.12). 

For the first term in
the right-hand side of (4.12) we observe: 
$$
P_+[P',Q_\lambda ]_+=(P[P',Q_\lambda ])_+-G^+(P)G^-([P',Q_\lambda ]).
$$
The s.g.o.\ term satisfies:
$$
G^+(P)G^-([P',Q_\lambda ])=\tfrac1\lambda G^+(P)G^-([P',P_1Q_\lambda
]),
$$
in view of (4.3). This shows (4.7) and (4.9), and (4.8) and (4.10)
follow by 
calculations such as:$$
G^-(P''Q_\lambda )=Jr^-P''(e^+r^++e^-JJr^-)Q_\lambda
e^+=G^-(P'')Q_{\lambda ,+}+\overline{P''}_+G^-(Q_\lambda ).\tag4.13
$$
Here $\overline{P''}=JP''J$ is likewise a $\psi $do satisfying the
transmission condition, and the calculation holds regardless of the
normal order of $P''$. We give details for the formulas in (4.10):
$$\aligned
G^-(P'P_1Q_\lambda )&=G^-(P'P_1)Q_{\lambda
,+}+(\overline{P'P_1})_+G^-(Q_\lambda ),\\
G^-(P_1Q_\lambda P')&=P_{1,+}G^-(Q_\lambda P')=P_{1,+}G^-(Q_\lambda
)P'_{+}+
P_{1,+}(\overline{Q_\lambda })_+G^-(P' ).
\endaligned$$
\qed
\enddemo

\proclaim{Lemma 4.2} Hypotheses as in Lemma {\rm
4.1}. 

$1^\circ$ The singular Green terms appearing in the primary formulas in
Lemma {\rm 4.1} are all of one of the forms$$
A''Q_{\lambda ,+},\quad 
A''G^{\pm}(Q_{\lambda }), \quad A''Q_{\lambda ,+}A''',
\quad A''G^{\pm}(Q_{\lambda })A''',
\tag4.14$$
(or the same expressions with $Q_\lambda $ replaced by $\overline
Q_\lambda $), where $A''=P''_++G''$ and $A'''=P'''_++G'''$ are of
normal order $0$ and class $0$. 

$2^\circ$ The singular Green terms appearing in the secondary formulas in
Lemma {\rm 4.1} (without the factor $\frac1\lambda $) are all of one
of the forms {\rm (4.14)} 
(or the same expressions with $Q_\lambda $ replaced by $\overline
Q_\lambda $), where $A''=P''_++G''$ has $P''$ of normal order $m$ and
$G''$ of class
$m$. The right factor 
$A'''=P'''_++G'''$ is of normal order $0$ and class $0$.  
The only resulting terms where $A''=P''_+$ of normal order $m$ occurs,
are of the form
$$
( PP_1)_+Q_{\lambda ,+}G'\text{ or }(PP_1)_+G^+(Q_\lambda
)G^-(P'),\tag4.15
$$ 
with an s.g.o.\ to the right.

\endproclaim 

\demo{Proof} For the terms in (4.5),  this is clear from the basic
rules of calculus, cf.\ e.g.\ \cite{G1}. For the terms in (4.6), $$
\gathered
GG'P_{1,+}Q_{\lambda
,+},\quad GP_{1,+}Q_{\lambda ,+}G', \quad P_+G'P_{1,+}Q_{\lambda
,+},\\
 P_+P_{1,+}Q_{\lambda ,+}G',\quad GP'_+P_{1,+}Q_{\lambda
,+},\quad GP_{1,+}Q_{\lambda ,+}P'_+,\endgathered
\tag 4.16$$
all the expressions except the fourth one have s.g.o.s of class $m$
to the left of $Q_{\lambda ,+}$, since, when $P_{1,+}$ is composed to the
left with an s.g.o.\ of class 0, we get an s.g.o.\ of class $m$. 
For the fourth expression, we observe that 
that since $P_1$ is a differential operator of order $m$,$$
P_+P_{1,+}=(PP_1)_++\sum_{0\le j\le {m-1}}K_j\gamma
_j,\tag4.17 
$$
where $PP_1$ is of normal order $m$ and $\sum_{0\le j\le {m-1}}K_j\gamma
_j$ is an s.g.o.\ of
class $m$; here the $K_j$ are Poisson operators of order $\sigma
+m-j$ and the $\gamma _j$ are the standard trace operators ($\gamma _ju=(D_n^ju)|_{x_n=0}$).
Thus $$
 P_+P_{1,+}Q_{\lambda ,+}G'=( PP_1)_+Q_{\lambda ,+}G'+G''Q_{\lambda ,+}G',
\tag4.18$$
where $PP_1$ has normal order $m$ and $G''$ has class $m$.

Now consider the terms in (4.7). The second and third term
are clearly of the asserted form with $A''$ of normal order and class
0. For the first term we use the
decomposition of $G^-([P',Q_\lambda ])$ given in (4.8) to reach this
conclusion.

Finally, consider the expressions in (4.9) with the additional
decomposition of a factor in the first term given in (4.10):
$$
\gather
G^+(P)[G^-(P'P_1)Q_{\lambda ,+}+(\overline {PP_1})_+G^-(Q_\lambda )
-P_{1,+}G^-(Q_\lambda
)P'_{+}-
P_{1,+}(\overline{Q_\lambda })_+G^-(P' )],
\\
G^+(P')\overline P_{1,+}G^-(Q_\lambda ),\quad P_+P_{1,+}G^+(Q_\lambda
)G^-(P')\tag4.19
\endgather
$$ All the
expressions have $P_1$ or $\overline P_1$ entering in compositions to
the left of a $\lambda $-dependent factor.
The first line and the first
expression in the second line lead to expressions with s.g.o.s
of class $m$ to the left. The last expression leads in view of (4.17) to
$$
P_+P_{1,+}G^+(Q_\lambda
)G^-(P')=   (PP_1)_+G^+(Q_\lambda)G^-(P')+G''G^+(Q_\lambda)G^-(P',)
$$
where $PP_1$ has normal order $m$ and $G''$ has class $m$.
\qed

\enddemo

We now investigate the normal traces.

\proclaim{Proposition 4.3} 
Let $G_\lambda $ be a parameter-dependent
singular Green operator of a form as
in Lemma {\rm 4.3}, and such that the sum of the
orders of $A''$ and 
$A'''$ is
$\varrho $; then $S_\lambda =\tr_n G_\lambda $ is a $\psi $do in the
parameter-dependent calculus of order $\varrho -m$ with
symbol $s(x',\xi ',\lambda
)\sim \sum_{j\ge 0}s_{\varrho -m-j}(x',\xi ',\lambda )$.

$1^\circ$ When $A''$ is of normal order and class $0$,
$S_\lambda $ is of regularity $\varrho -\frac14$, and the symbol 
satisfies estimates (where $\mu =|\lambda |^{\frac1m})$:
$$
|\partial _{x',\xi '}^{\beta ,\alpha } [s(x',\xi ',\lambda
)-\sum_{j<J}s_{\varrho -m-j}(x',\xi ',\lambda )]|
\leg \ang{\xi '}^{\varrho -\frac14 -|\alpha |-J}\ang{\xi ',\mu
}^{-m+\frac14},
\tag4.20$$
for all indices.

$2^\circ$ When $A''$ is of normal order $m$ and class $m$,
$S_\lambda $ is of regularity $\varrho -m+\frac14 $, the symbol 
satisfying estimates :
$$
|\partial _{x',\xi '}^{\beta ,\alpha } [s(x',\xi ',\lambda
)-\sum_{j<J}s_{\varrho -m-j}(x',\xi ',\lambda )]|
\leg \ang{\xi '}^{\varrho -m+\frac14-|\alpha |-J}\ang{\xi ',\mu
}^{-\frac14},
\tag4.21
$$
for all indices.

\endproclaim

\demo{Proof} 
Note that $G_\lambda $ is in all cases of class 0, since $Q_{\lambda
}$ is of order $-m$ and $A'''$ is of normal order and class 0.

Consider first the case where there is no factor $A'''$ (or when
$A'''=I$); here we get the results
fairly easily. 
Then $A''$ is of order $\varrho $. If the class of $G''$ is
0, $A''$ enters in the parameter-dependent calculus as the sum of a $\psi
$do and an s.g.o.\, both of
order $\varrho $ and regularity $\varrho $, so when it is composed
with $Q_{\lambda ,+}$ or $G^{\pm}(Q_\lambda )$ of order $-m$ and regularity
$+\infty $ we get an operator of order $\varrho -m$ and regularity
$\varrho $, in view of \cite{G1, Th.\ 2.7.7, Cor.\ 2.7.8}. By Lemma 3.4,
$\tr_n$ of it is a $\psi $do on $X'$ of order 
$\varrho -m$ and regularity $\varrho -\frac14$.

When $\varrho \le \frac14$, the estimates (4.20) hold automatically
(are standard symbol estimates),
since the power of $\ang{\xi ',\mu }$ in the parenthesis in (3.6) can
be left out when $\nu \le 0$. For larger $\varrho $, we compose to
the left with $\Lambda ^\varrho \Lambda ^{-\varrho }$, where $\Lambda
^t=\op '(\ang{\xi '}^t)$; it is accounted for in \cite{G1, Sect.\ 2.8}
that this defines operators within the calculus (not only for $\psi
$do's on $\Bbb R^{n-1}$ but also for s.g.o.s). The preceding
considerations now apply to the expression composed to the left with
$\Lambda ^{-\varrho }$, which satisfies estimates with an extra factor
$\ang{\xi '}^{-\varrho }$. The resulting operator will satisfy
(4.20) with $\varrho $ replaced by 0, and when we recompose with
$\Lambda ^\varrho $ to the left, it is easily checked from the
composition rules that we obtain an operator satisfying (4.20).

When $A''=G''$ of class $m$, its regularity is only $\varrho -m+\frac12$
(by \cite{G1, (2.3.55)}); then the composed operator has regularity
$\varrho -m+\frac12$, and $\tr_n$ of it has regularity $\varrho
-m+\frac14 $ by Lemma 3.4. 
(The central fact here is that $G''=G_1+\sum_{l<m}K_l\gamma _l$, such
that the term with the
weakest decrease in $\lambda $ will be $K_{m-1}\gamma
_{m-1}Q_{\lambda ,+}$, where $\gamma _{m-1}Q_{\lambda ,+}$ is a trace
operator whose symbol norm is $O(\ang{\xi ',\mu }^{-\frac12}$).)
If $\varrho -m<-\frac14$, the estimates
(4.21) are automatically satisfied, otherwise we obtain them by
pulling out a factor $\Lambda ^t\Lambda ^{-t}$ for a large $t$ as above.

Next, we consider the case where $A'''$ is nontrivial. Here we have
to make some extra 
efforts, both since regularity numbers in compositions are not in
general additive, and since we have to deal with some inconvenient
terms (4.15). There are now
three factors, with the $\lambda $-dependent factor in the middle.

 Let $A''$ and
$A'''$ have 
orders $\varrho _1$ and $\varrho _2$, so that $\varrho =\varrho
_1+\varrho _2$.
Invoking the trick of composing to the left with $\Lambda ^{\varrho
_1}\Lambda ^{-\varrho _1}$ if $\varrho _1>0$ and to the right with
$\Lambda ^{\varrho _2}\Lambda ^{-\varrho _2}$ if $\varrho _2>0$, we
can assume that $\varrho _1,\varrho _2\le 0$. Assume to begin with
also that $q$ is independent of $x_n$.

The normal trace of $G_\lambda $ is found by applying (3.15)
to its symbol-kernel \linebreak$g(x',\xi ',x_n,y_n,\lambda )$.
We recall from \cite{G1} the notation $g(x',\xi ',D_n,\lambda )$
(or just $g(D_n)$) for
the operator on $\rp$ defined for each $(x',\xi ',\lambda )$ by
applying the s.g.o.\ definition in one variable $x_n$; we use again
the notation $\circ_n$ for the composition
of such one-dimensional operators. For operators on $\Bbb R_+$ of
normal order and  class 0, $\tr_n$ is the usual trace, so there is a
certain commutativity, namely e.g.\ 
$$
\aligned
\tr_n\bigl(g(D_n )\circ_n g'(D_n)\bigr)&=
\tr_n\bigl(g'(D_n )\circ_n g(D_n )\bigr),\\
\tr_n\bigl(g(D_n )\circ_n p(D_n)_+\bigr)&=
\tr_n\bigl(p(D_n )_+\circ_n g(D_n )\bigr);
\endaligned\tag 4.22$$
the s.g.o.s are
smoothing.
We shall use this to reduce the most difficult estimates for three
components to cases of two components with better properties. Consider e.g.\ a
composition $A''Q_{\lambda ,+}A'''$. Here$$
\tr_n(a''\circ_n q_+\circ_n a''')=\tr_n(a'''\circ_n a''\circ_n
q_+),\tag 4.23
$$ 
since $a''\circ_n q_+$  and $ a'''$ are both of normal order
and class 0. In the compositions coming from the primary cases in
Lemma 4.1, $a'''\circ_n a''$ will be of normal order and class 0.
In the compositions coming from the secondary cases in
Lemma 4.1, $a'''\circ_n a''$ will be {\it a singular Green operator
of class} $m$; this is clear if $a''$ is such one, and if
$a''=p''_+$ is of normal order $m$, $a'''$ is necessarily an s.g.o.\ of
class 0 according to Lemma 4.2, so the composite is an s.g.o.\ of
class $m$.  The important 
fact is that we get rid of contributions of the form $(pp_1)_+q_+g'$,
where a direct attack need not give estimates with a decrease in
$\lambda $ since 
$$
\sup_{\xi _n}|p_{1,m}(x,\xi )(p_{1,m}(x,\xi )-\lambda )^{-1}|=1.\tag 4.24
$$
 Now the results from the beginning of the proof for compositions
$A''Q_{\lambda ,+}$ can be applied.
When $a'''\circ _n a''$ is of normal order and class 0, this gives a
symbol of order $\varrho -m$ and regularity $\varrho -\frac14$, and
when $a'''\circ _n a''$ is an s.g.o.\ of class $m$, we get a
symbol of order $\varrho -m$ and regularity $\varrho -m+\frac14$;
since $\varrho \le 0$, the estimates in (4.20)--(4.21) are automatic.

The commutation is only allowed on the one-dimensional level. To find
the full composition of $a''$, $q_+$ and $a'''$, we note that$$
 \tr_n
(a''\circ q_+\circ a''')
\sim\sum_{\alpha , \beta \in \Bbb N^{n-1}}\frac{(-i)^{\alpha +\beta
}}{\alpha !\beta !}\tr_n(\partial _{\xi '}^{\alpha } a''\circ_n 
\partial _{x'}^{\alpha}\partial _{\xi '}^\beta q_+
\circ_n \partial _{x'}^{\alpha }\partial _{x'}^\beta  a'''), 
\tag4.25$$
and perform the above commutation idea for each term, to find the
desired symbol information. 

Concerning
remainders, an analysis shows that it is only the part of normal
order $m$ of $PP_1$, giving a term  $D'D_n^m
Q_{\lambda ,+}G'''$, that needs special treatment; for the part
$P_{m-1}$ of 
$PP_1$ of normal order $\le m-1$ one can appeal to the estimate
$$
\sup_{\xi _n}|p_{m-1}(p_{1,m}-\lambda )^{-1}|\leg \ang{\xi
'}^{\sigma +1}\ang{\xi ',\mu }^{-1}.\tag4.26
$$
In the usual remainder term (as in e.g.\ \cite{GSc2, pf.\ of Prop.\
3.8}) in the calculation of the composition
inside 
$\tr_n((D_n^m
Q_{\lambda })_+G''')$, one can then
perform a commutation (4.22) inside the integral w.r.t.\  $h$.

 If $q$ depends on $x_n$, it must be Taylor expanded
in $x_n$ and each term treated individually; here one uses that in
the terms with $k$'th powers of
$x_n$, $k\ge 1$, the symbols coming from $q$ are $O(\lambda ^{-2})$ and the
order of the s.g.o.s are lowered by $k$.

There is a similar analysis when $Q_{\lambda ,+}$ is replaced by
$\overline Q_{\lambda ,+}$, $G^{\pm}(Q_{\lambda })$ or
$G^{\pm}(\overline Q_{\lambda })$.
\qed

\enddemo  

It may be remarked that the fraction $\frac14$ comes in because of
the general application of Lemma 3.4. Particular efforts applied to
the individual compositions may give an
improvement to $\frac 
12$ in  (4.21) --- and an analysis extending that
of \cite{GSc1} would give further improvements, cf.\ Remark 3.3. But
the gain of $\frac14$ is sufficient for the present purposes.  

We can now conclude:

\proclaim{Theorem 4.4} Let $A=P_++G$ of order $\sigma $ and normal
order and class $0$, let $A'=P'_++G'$ of order $\sigma '$ and normal
order and class $0$, and let $P_1$ be an auxiliary elliptic
differential operator  
of order $m>\sigma +\sigma '+n$, with no eigenvalues of the
principal symbol on $\crm$ (so that $Q_\lambda =(P_1- \lambda )^{-1}$
is defined for large $\lambda $ in a sector $V$ around $\rmi$). We
assume that $P$ and $P'$ are zero if $\sigma $ or 
$\sigma '$ is noninteger.

Let $ {\Cal S}_\lambda
=\tr_n {\Cal G}_\lambda $ with symbol $\frak s(x',\xi ',\lambda )$,
where $\Cal G_\lambda $ is the singular Green part of
\linebreak$A[A',Q_{\lambda ,+}]$.
Then
$ {\Cal S}_\lambda $ is a family of $\psi
$do's on $X'$ with the
properties:

$1^\circ$ $ {\Cal S}_\lambda $ is of order $\sigma +\sigma '-m$ and regularity
$\sigma +\sigma '-\frac14$, the symbol satisfying:$$
|\partial _{x',\xi '}^{\beta ,\alpha } [\frak s(x',\xi ',\lambda
)-\sum_{j<J}\frak s_{\sigma +\sigma ' -m-j}(x',\xi ',\lambda )]|
\leg \ang{\xi '}^{\sigma +\sigma ' -\frac14 -|\alpha |-J}\ang{\xi ',\mu
}^{-m+\frac14},
\tag4.27$$
on the rays in $V$ (with $-\lambda =\mu ^me^{i\theta }$, $\mu >0$), for all $\alpha
,\beta ,J$.

$2^\circ$ $\lambda  {\Cal S}_\lambda $ is of order $\sigma
+\sigma '$ and regularity 
$\sigma +\sigma '+\frac14$, and  for all $\alpha ,\beta ,J$,
$$
|\partial _{x',\xi '}^{\beta ,\alpha } [\frak s(x',\xi ',\lambda
)-\sum_{j<J}\frak s_{\sigma +\sigma ' -m-j}(x',\xi ',\lambda )]|
\leg \ang{\xi '}^{\sigma +\sigma ' +\frac14-|\alpha |-J}\ang{\xi ',\mu
}^{-\frac14}\mu ^{-m}.
\tag4.28
$$

\endproclaim 

\demo{Proof} 
This follows immediately from Proposition 4.3 in view of the
description of $ {\Cal G}_\lambda $ given in Lemmas 4.1 and 4.2.\qed
\enddemo 

We can then establish trace expansions. Here we first consider the
case where $\sigma $ and $\sigma '$ are integers. 

\proclaim{Theorem 4.5} Assumptions as in Theorem {\rm 4.4}, with
$\sigma $ and $\sigma '\in\Bbb Z$.

There is a trace expansion
$$
\Tr([A,A']Q_{\lambda ,+})=\sum_{0\le j \le n+\sigma +\sigma
'}c_{j}\,(-\lambda ) ^{\frac{n+\sigma +\sigma '-j}m-1}+ O(\lambda
^{-1-\frac1{4m}}),\tag4.29
$$
so that
$$
C_0([A,A'], P_{1,+})=c_{n+\sigma +\sigma '}\tag4.30
$$
(taken equal to $0$
if $n+\sigma +\sigma '<0$) is well-defined.

The
symbol $s(x',\xi ')$ deduced from the symbol $\frak s(x',\xi
',\lambda )$ of $ {\Cal S}_\lambda $ by $$
s(x',\xi ')=\tfrac i{2\pi }\int_{\Cal C''}\log\lambda \, \frak
s(x',\xi ',\lambda )\, d\lambda \tag4.31
$$
(with $\Cal C''$ a curve in $\Bbb C\setminus\crm$ encircling the
sectorial set containing the eigenvalues of $p_{1,m}(x,\xi )$ for
$x\in X $, $|\xi '|\ge 1$), 
is a classical $\psi $do symbol of order $\sigma +\sigma '$, defining
a $\psi $do $S$ such that$$
C_0([A,A'],P_{1,+})=-\tfrac1m\operatorname{res}((P[P',\log P_1
])_+)-\tfrac1m\operatorname{res}(S).\tag 4.32
$$
\endproclaim 

\demo{Proof} 
For $\Cal R_\lambda =P[P',Q_\lambda ]$ we have a diagonal kernel
expansion as in (2.22)ff.\ with coefficients $\tilde r_j(x)$. Integrating over
the coordinate patches 
intersected with $\crnp$, we find that
$$
\aligned
\Tr\Cal R _\lambda 
&=
\sum_{ j<\sigma +\sigma '+m+n }   \tilde r_{ j,+}(-\lambda )
^{\frac{n+\sigma +\sigma ' -j}m  -1}+ 
O(|\lambda |^{-2+\varepsilon }),\text{ where}\\
\tilde r_{j,+}&=\int_{\rnp}\tr \tilde r_j(x),\quad \tilde r_j(x)=  \int_{\Bbb R^n}r^h_{-m+\sigma +\sigma '-j}(x,\xi
,-1)\,\d\xi .
\endaligned \tag4.33
$$
The calculations around Theorem 2.3 apply to this situation,
showing that the coefficient of $(-\lambda )^{-1}$ identifies with
the residue:
$$
\tilde r_{n+\sigma +\sigma ',+}=-\tfrac1m
\operatorname{res}\bigl((P[P',\log P_1])_+\bigr).\tag4.34
$$

For $ {\Cal S}_\lambda $, the information that it
is of order $\sigma +\sigma '-m$ and regularity
$\sigma +\sigma '-\frac14$ leads by Lemma 3.5 to a trace expansion
$$
\aligned
\Tr_{\Bbb R^{n-1}} {\Cal S}_\lambda  
&=\sum_{1\le j < n+\sigma +\sigma '-\frac14} \frak s_{j}\,(-\lambda )
^{\frac{n+\sigma +\sigma '-j}m-1}+O(\lambda
^{-1+\frac1{4m}}), \\
\frak s_j&=\tilde {\frak s}_{j-1}=\int\tr \tilde{\frak s}_{j-1}(x')\,dx',\\
\tilde {\frak s}_{l}(x')&=
\int_{\Bbb R^{n-1}}{\frak s}^h_{\sigma +\sigma '-m-l}(x',\xi ',-1)\,\d\xi',
\endaligned
\tag4.35
$$ 
which just misses having a precise term $c(-\lambda
)^{-1}$.
But we can improve the expansion by using the additional information
we have on the symbol in Theorem 4.4.
In fact, for $j=\sigma +\sigma '+n$,  $l=\sigma +\sigma '+n-1$, we have
a term
(taken equal to 0 if
$\sigma +\sigma '\le -n$) satisfying
$$\aligned
|\frak s_{-m-n+1}(x',\xi ',\lambda )|&\leg \ang{\xi '}^{-\frac14-n+1}
\ang{\xi ',\mu
}^{-m+\frac14},\\
|\frak s_{-m-n+1}(x',\xi ',\lambda )|&\leg \ang{\xi '}^{\frac14-n+1}
\ang{\xi ',\mu
}^{-\frac14}\mu ^{-m},
\endaligned\tag4.36$$
and the remainder $\frak s'=\frak s-\sum_{l< \sigma +\sigma '+n}\frak
s_{\sigma +\sigma '-m-l+1}$ after this term satisfies
$$\aligned
|\frak s'|&\leg \ang{\xi '}^{-\frac14-n}\ang{\xi ',\mu
}^{-m+\frac14},\\
|\frak s'|
&\leg \ang{\xi '}^{\frac14-n}\ang{\xi ',\mu
}^{-\frac14}\mu ^{-m}.
\endaligned\tag4.37$$
From (4.36) follows as in \cite{G1, Lemma 2.1.9} that
$$\aligned
|\frak s^h_{-m-n+1}(x',\xi ',\lambda )|&\leg |\xi '|^{-\frac14-n+1}
|\xi ',\mu|^{-m+\frac14},\\
|\frak s^h_{-m-n+1}(x',\xi ',\lambda )|&\leg |\xi '|^{\frac14-n+1}
|\xi ',\mu
|^{-\frac14}\mu ^{-m},
\endaligned\tag4.38$$
so $\frak s^h_{-m-n+1}$ is integrable at $\xi '=0$ (besides being so
for $|\xi '|\to \infty $) when $\lambda  \ne 0$. Then 
$$
\aligned
\Tr(\operatorname{OP}'(\frak s^h_{-m-n+1}))&=\frak s_{\sigma +\sigma
'+n}\,(-\lambda )^{-1},\text{ with}\\
\frak s_{\sigma +\sigma '+n}&=\tilde {\frak s}_{\sigma +\sigma
'+n-1}=\int\tr \tilde{\frak s}_{\sigma +\sigma '+n-1}(x')\,dx',\\
\tilde {\frak s}_{\sigma +\sigma '+n-1}(x')&=
\int_{\Bbb R^{n-1}}{\frak s}^h_{-m-n+1}(x',\xi ',-1)\,\d\xi',
\endaligned\tag4.39
$$
as in (4.35). This gives the needed extra term, but we also have to show
that remainders do not interfere. 
(4.28) shows that $|\frak s'| \leg \ang{\xi
'}^{\frac14-n}\mu ^{-m-\frac14}$, which integrates in $(n-1)$-space to
give an estimate by $\mu ^{-m-\frac14}$. The difference $\frak
s^h_{-m-n+1}-\frak s_{-m-n+1}$ is $O(\mu ^{-m-\frac14})$ on its
support contained in $\{|\xi |\le 1\}$, so it likewise integrates to
an $O(\mu ^{-m-\frac14})$ term. This also holds for the preceding
terms, the differences
$\frak
s^h_{\sigma +\sigma '-m-l}-\frak s_{\sigma +\sigma '-m-l}$ with $
l<\sigma +\sigma '+n-1$. Then we can finally
conclude (4.29). 

We shall now show that the integral in (4.31) is
well-defined so that the symbol properties can be checked directly.
Again we use the estimates in Theorem 4.4.
Note that (4.27) gives too little
decrease in $\lambda  $ to allow the integration (4.31), whereas
(4.28) gives enough decrease in $\lambda  $, but much  less in
$\xi '$. 
Using that $\frak s$, its terms
and remainders are $O(\lambda ^{-1-\frac14 })$, we can insert $\frak
s$ in (4.31) in order to obtain $s(x',\xi ')\sim 
\sum_{j\ge 0}s_{\sigma +\sigma '-j}(x',\xi ')$. Here
$s_{\sigma +\sigma '-j}$ is homogeneous of degree $\sigma +\sigma
'-j$ in $\xi '$ for $|\xi '|\ge 1$, in view of the following calculation
with  $t\ge 1$, $\varrho =t^{-m}\lambda $:
$$\aligned
s_{\sigma +\sigma '-j}&(x',t\xi ')=\tfrac i{2\pi }\int_{\Cal
C''}\frak s_{\sigma +\sigma '-m-j}(x',t\xi ',\lambda
)\log\lambda \,d\lambda \\
&=\tfrac i{2\pi }\int_{\Cal C''}t^{\sigma +\sigma '-m-j}\frak
s_{\sigma -\sigma '-m-j}(x',\xi ',\varrho 
)(\log\varrho +m\log t)t^m \,d\varrho  \\
&=t^{\sigma +\sigma '-j}\tfrac i{2\pi }\int_{\Cal C''}\frak s_{\sigma
+\sigma '-m-j}(x',\xi ',\varrho 
)\log\varrho  \,d\varrho =t^{\sigma+\sigma '-j}s_{\sigma +\sigma '-j}(x',\xi '),  
\endaligned\tag 4.40$$ 
where we have used that $\tfrac i{2\pi }\int_{\Cal C''}\frak
s_{\sigma +\sigma '-m-j}(x',\xi
',\varrho  
) \,d\varrho =0$ since the integrand is holomorphic on the region to
the left of $\Cal C''$ and
$O(\lambda ^{-\frac54 })$ for $\lambda \to\infty $ there. 

Remainders satisfy$$
|s(x',\xi ')-\sum_{ j<J}s_{\sigma +\sigma '-j}(x',\xi ')|\leg
\ang{\xi '}^{\sigma +\sigma '+\frac14-J}\tag4.41
$$
for  all $J$, in view of (4.28).
Using the exact 
terms for $ j<J'=J+1$ and the remainder estimate
(4.41) with $J$ replaced by $J'$, we can improve (4.41) to
 $$
|s(x',\xi ')-\sum_{j<J}s_{\sigma +\sigma '-j}(x',\xi ')|\leg
\ang{\xi '}^{\sigma +\sigma '-J},\tag4.42
$$
which is the appropriate estimate for showing that $s$ is
polyhomogeneous of order $\sigma +\sigma '$. Estimates of derivatives are included
in a similar way.

So now $s$ is well-defined as a classical symbol of order $\sigma
+\sigma '$; it
defines the operator $S$ with the residue $$
\operatorname{res}S=\int_{\Bbb R^{n-1}}\int_{|\xi '|=1}\tr s_{1-n}(x',\xi ')\,\d
S(\xi ')dx'.\tag 4.43
$$
From the fact that 
$$
s_{1-n}(x',\xi ')=\tfrac i{2\pi }\int_{\Cal C''}\log\lambda\, \frak
s^h_{-m-n+1}(x',\xi ',\lambda )\, d\lambda 
$$
for $|\xi '|\ge 1$, it is found by use of Lemma 1.2 and Lemma 1.3
for dimension $n-1$,  that 
$$
-\tfrac1m\operatorname{res}S=\int_{\Bbb R^{n-1}}\int_{\Bbb
R^{n-1}}\frak s_{-m-n+1}(x',\xi ',-1)\,\d\xi 'dx'=\frak s_{\sigma +\sigma '+n}.\tag4.44
$$
Collecting the residues and contributions to $C_0([A,A'], P_{1,+})$ from
(4.34) and (4.44), we find (4.32).
\qed
\enddemo 

Noninteger orders are included as follows:

\proclaim{Theorem 4.6} 
Assumptions as in Theorem {\rm 4.4}, with
$\sigma $ and $\sigma '\in \Bbb R$ and $P=P'=0$.

There is a trace expansion
$$
\Tr([A,A']Q_{\lambda ,+})=\sum_{0\le j < n+\sigma +\sigma '+\frac14
}c_{j}\,(-\lambda ) ^{\frac{n+\sigma +\sigma '-j}m-1}+ O(\lambda
^{-1-\frac1{4m}\;(+\varepsilon )}),\tag4.45
$$
where $\varepsilon =0$ if $\sigma +\sigma '+\frac14\notin\Bbb Z$.
Define
$
C_0([A,A'], P_{1,+})=c_{n+\sigma +\sigma '}
$ if $n+\sigma +\sigma '\in\Bbb N$, $C_0([A,A'], P_{1,+})=0$
otherwise.  Then
defining $S$ as in Theorem {\rm 4.8}, we have that
$$
C_0([A,A'],P_{1,+})=-\tfrac1m\operatorname{res}(S).\tag 4.46
$$
\endproclaim 

\demo{Proof} There is no $\psi $do term in this case. For the s.g.o.\
term $ {\Cal G}_\lambda $ we proceed as in the preceding
proof. It goes 
over verbatim if $\sigma +\sigma '\in\Bbb Z$, whereas one has to
modify the indexations when $\sigma +\sigma '\notin\Bbb Z$. Actually,
that is a case where there will be no nontrivial term  $c(-\lambda
)^{-1}$, and all one has to check is remainder estimates. Since $S$ is
of noninteger order then, $\operatorname{res}S$ is also zero.\qed
\enddemo

\subhead 5. Extension of the res of log formula to pseudodifferential
boundary problems \endsubhead

With these techniques at hand, we shall also investigate possible
extensions of the res of log formula (1.2) to
realizations of elliptic pseudodifferential boundary problems.
Consider a normal elliptic
realization $(P_++G)_T$, as defined in \cite{G1, Section 3.3}. Here
$P$ is a classical $\psi $do in $E$ of integer order $m>0$ satisfying the
transmission condition at $X'$, $G$ is a singular Green operator in $E$
of order and class $m$, and $T=\{T_0,\dots,T_{m-1}\}$ is a normal
trace operator with 
entries $T_k$ of order and class $k$ going from $E$ to $F_k$,
all polyhomogeneous. $E$ and the $F_k$  are hermitian $C^\infty $
vector bundles 
over $X$ resp.\ $X'$.
We assume that the conditions for uniform
parameter-ellipticity in \cite{G1, Def.\ 3.3.1} are satisfied on the
rays in a sector $V$ around $\rmi$.

The resolvent$$
((P_++G)_T-\lambda )^{-1}=R_\lambda =Q_{\lambda ,+}+G_\lambda \tag5.1
$$
was constructed in \cite{G1, Sect.\ 3.3} and shown to belong to the
parameter-dependent calculus set up in the book. Complex powers
$((P_++G)_T)^z$ were described to some extent in \cite{G1, Sect.\ 4.4},
just for $\operatorname{Re}z<0$, where it was shown that their
singular Green part has some, but not all of, the symbol estimates of
standard s.g.o.s. The logarithm of $(P_++G)_T$ has not, to our
knowledge, been discussed anywhere. 

Since the complex powers were only considered for
$\operatorname{Re}z<0$, we cannot draw conclusions about a derivative
at $z=0$, but one can try a formula as in (1.14); it generally leads
to an operator outside the Boutet de Monvel calculus. Rather than 
going into a deeper analysis of such operators and the possibility of
defining residues on them, we shall show a generalization of (1.2)
where a residue of the logarithm of the $\psi $do part does enter,
and the s.g.o.\ part is reduced to the residue of a classical $\psi
$do on $X'$; the ``nice part'' of the log contribution from
$G_\lambda $.

It is shown in \cite{G1, Th.\
3.3.5, 3.3.10} that when $m>n$, the resolvent
has a trace expansion with at least $n+1$ exact terms:
$$
\Tr R_\lambda =\sum_{0\le j\le n}c_j(-\lambda
)^{\frac{n-j}m-1}+O(\lambda ^{-1-\frac1{4m}}),\tag5.2 
$$
valid for $\lambda \to\infty $ in the sector of
parameter-ellipticity. (If the regularity is greater than $1$,
there will be more terms in the expansion.) The coefficients $c_j$ are
defined by 
integration of the strictly homogeneous terms in the symbols of
$Q_{\lambda ,+} $ and $G_\lambda $; in particular, the coefficient of
$(-\lambda )^{-1}$,
$$C_0(I,(P_++G)_T)=c_n\tag5.3
$$
is defined from the term of order $-m-n$ in the symbol of $Q_{\lambda
,+}$ and the term of order $-m+1-n$ in the symbol of $G_\lambda $ (in
local coordinates).
As usual, $Q_\lambda $ is the inverse of $P-\lambda $, defined on a
larger compact 
$n$-dimensional manifold $\widetilde X$ in which $X$ is smoothly
imbedded.

In the following, we work in a localization to
 $\Bbb R^n$ (with $X$ carried over to
subsets of $\crnp$), as in the preceding sections. Let $Q_\lambda $,
$G_\lambda $ and 
$S_\lambda =\tr_nG_\lambda $ have
symbols $q$,  $ g$ and $s=\tr_ng$, respectively, with expansions e.g.$$
\aligned
q(x,\xi ,\lambda )&\sim\sum_{j\ge 0}q_{-m-j}(x,\xi ,\lambda ),\\ 
s(x',\xi ',\lambda )&\sim\sum_{j\ge 0}s_{-m-j}(x',\xi ',\lambda
).
\endaligned\tag5.4
$$Then
$$
\aligned
c_n&=c^P_{n,+}+c^G_n,\text{ with}\\
c^P_{n,+}&=\int_{\rnp}\int_{\Bbb R^n}\tr q^h_{-m-n}(x,\xi ,-1)\,\d\xi dx,\\
c^G_n&=\int_{\Bbb R^{n-1}}\int_{\Bbb R^{n-1}}\tr s^h_{-m+1-n}(x',\xi ',-1)\,\d\xi' dx'.
\endaligned\tag5.5
$$

Consider the elliptic system $\{P_++G,T\}$ defining the
operator $(P_++G)_T$ we are interested in. The order is $m$, the
regularity of $P$ is $m$, and the
regularity $\nu $ of the full system is an integer or half-integer
lying in the interval 
$[\frac12,m]$ (cf.\ \cite{G1, (3.3.11)}) --- unless the operators are
purely differential, in 
which case the regularity is $+\infty $ (any $\nu \in\Bbb R$ works then).
As shown in \cite{G1, Th.\ 3.3.2}, $Q_\lambda $ is of order $-m$ and
regularity $m$,  and
$G_\lambda $ is of order $-m$, class $0$ and regularity $\nu $ 
(the regularities being replaced by $+\infty $ in the differential
operator case).

With reference to the lemmas in Section 3 here,
the proof of (5.2) in \cite{G1, Sect.\ 3.3} consists of 
applying Lemma 3.2 to the pseudodifferential part $Q_{\lambda ,+}$
to get pointwise expansions of the diagonal kernel of $Q_\lambda $
and integrate 
these over $\rnp$,  applying Lemma 3.5 to the normal trace of
the s.g.o.\ part $G_\lambda $
to get pointwise expansions of the diagonal kernel and integrate
these over $\Bbb R^{n-1}$ (contributions from interior patches are
smoothing and $O(\lambda ^{-1-\frac1{4m}})$), and adding the
expansions.
  
\comment
Let us recall the arguments for (5.2) given in \cite{G1, Sect.\
3.3}, with reference to the lemmas in Section 3 here: Lemma 3.2
provides us with a diagonal kernel expansion, and hence 
by integration in $x\in\rnp$ a trace expansion of $Q_{\lambda ,+}$:
$$
\aligned
K(Q_\lambda ,x,x,\mu ) &=\sum_{0\le j< n+m}c^P_j(x)\mu
^{n-j-m}+O(\mu  ^{-2m+\varepsilon }),\\ 
\Tr Q_{\lambda ,+} &=\sum_{0\le j< n+m}c^P_{j,+}\mu 
^{n-j-m}+O(\mu  ^{-2m+\varepsilon }).
\endaligned\tag5.6
$$
Lemma 3.4 shows that $S_\lambda =\tr_n G_\lambda $ is of regularity
$\nu -\frac14$, and Lemma 3.5 gives the expansions of kernel and trace:
$$
\aligned
K(S_\lambda ,x',x',\mu ) &=\sum_{0\le l\le n-1+\nu -\frac14}\tilde s_l(x')\mu
^{n-1-l-m}+O(\mu  ^{-m-\nu +\frac14 }),\\ 
\Tr _{\rnp} G_\lambda =\Tr_{\Bbb R^{n-1}}S_\lambda  &=\sum_{0\le l\le
n-1+\nu -\frac14}\tilde s_l\mu
^{n-1-l-m}+O(\mu  ^{-m-\nu +\frac14 }).
\endaligned\tag5.7 
$$
The latter can also be written, when we replace $l$ by $j-1$ and
rename the coefficients,$$
\Tr _{\rnp} G_\lambda =\sum_{1\le j\le
n+\nu -\frac14}c^G_j\mu
^{n-j-m}+O(\mu  ^{-m-\nu +\frac14 }).\tag5.8
$$
Adding this to the last line in (5.6) and replacing $\mu $ by
$(-\lambda )^{\frac1m}$, we find (5.2) with (5.3), (5.5), valid for
$\lambda \in\rmi$. There are similar expansions on the other rays in
$V$, and they must have the same coefficients as here, in view of the
holomorphy.
\endcomment

Now we want to relate the coefficients $c^P_{n,+}$ and $c^G_n$ to
residues. $c^P_{n,+}$ is immediately understood on the basis of Theorem
1.3 (integrating the pointwise version over $\rnp$). For $c^G_n$, we have
the following lemma.

To explain the curve $\Cal C''$ used there, we recall from \cite{G1}
that the ellipticity hypothesis assures that the strictly homogeneous
principal symbol $p^h_m(x,\xi )-\lambda $ and principal boundary
symbol operator $\{p^h(x',0,\xi ',D_n)+g^h(x',\xi ',D_n)-\lambda ,
t^h(x',\xi ',D_n)\}$ are invertible for $\lambda $ in a sector around
$\rmi$, $(\xi ',\lambda )\ne 0$, such that the resolvent exists in a
keyhole region $V_{r,\varepsilon }$ (1.18) except at finitely many
points. By a small rotation, we can assure that no eigenvalues are on
$\rmi$. As $\Cal C''$ we take a curve in $\Bbb C\setminus \crm$
around $\complement V_{r,\varepsilon }$ and  the spectrum except
possibly 0; it can be the boundary of 
$V_{r',\varepsilon '}$ with suitably small $r'$ and $\varepsilon '$.

\proclaim{Lemma 5.1}
Define from $s$ and $S_\lambda $ the reduced symbol $s'$ and the corresponding operator
$S'_\lambda $:$$\aligned
s'(x',\xi ',\lambda )&=s(x',\xi ',\lambda ) - s_{-m}(x',\xi ',\lambda
),\\
S'_\lambda &=\operatorname{OP}'(s'(x',\xi ',\lambda )),
\endaligned\tag5.6$$
and set$$
\aligned
B&=\tfrac i{2\pi }\int_{\Cal C''}S'_\lambda \log\lambda \,d\lambda ,\\
b_{-j}(x',\xi ')&=\tfrac i{2\pi }\int_{\Cal C''}s_{-m-j}(x',\xi ',\lambda
)\log\lambda \,d\lambda \text{ for }j\ge 1, |\xi '|\ge 1.
\endaligned\tag5.7$$
Then $B$ is a classical $\psi $do on $\Bbb R^{n-1}$ of order $-1$ with symbol
$b\sim\sum_{j\ge 1}b_{-j}$.
\endproclaim 

\demo{Proof} 
Since $s'$ is of order $-m-1$ and regularity $\nu -\frac54$,
we have that
$$s' \text{ is }O((\ang{\xi '}^{\nu -\frac 54}+\ang{\xi ',\mu }^{\nu
-\frac54})\ang{\xi 
',\mu }^{-m-1-\nu +\frac54}),
$$ hence falls off like $\lambda 
$ to the power $\max\{-1-\frac1m,-1-\frac{\nu -\frac14}m\}$, so the symbol
multiplied by $\log \lambda $ is $O(\lambda ^{-1-\delta })$ with
a $\delta >0$. There are similar estimates for derivatives. Then $B$
is defined as a bounded operator in $L_2$, and its  symbol
terms $b_{-j}$ are found by integration of the terms in $s'$ as
stated.
To see that $b_{-j}$ is homogeneous of degree $-j$ in $\xi '$ for $|\xi
'|\ge 1$, we write for $t\ge 1$, with $\varrho =t^{-m}\lambda $:
$$\aligned
b_{-j}(x',t\xi ')&=\tfrac i{2\pi }\int_{\Cal C''}s_{-m-j}(x',t\xi ',\lambda
)\log\lambda \,d\lambda \\
&=\tfrac i{2\pi }\int_{\Cal C''}t^{-m-j}s_{-m-j}(x',\xi ',\varrho 
)(\log\varrho +m\log t)t^m \,d\varrho  \\
&=t^{-j}\tfrac i{2\pi }\int_{\Cal C''}s_{-m-j}(x',\xi ',\varrho 
)\log\varrho  \,d\varrho =t^{-j}b_{-j}(x',\xi '),  
\endaligned\tag 5.8$$
where the term with $m\log t$ drops out as in (4.40).
\comment
[er set f\o{}r, forkortes]
 
where we have used that $\tfrac i{2\pi }\int_{\Cal C''}s_{-m-j}(x',\xi
',\varrho  
) \,d\varrho =0$ since the integrand is holomorphic on the region to
the left of $\Cal C''$ and
$O(\lambda ^{-1-\delta })$ for $\lambda \to\infty $ there.
\endcomment  
Derivatives in $x'$ and $\xi '$ and remainders are easily checked.
\qed
\enddemo 

$B$ can in a sense be considered as the ``nice $\psi $do part'' of
the logarithmic 
contribution from the normal trace of the singular Green term
$G_\lambda $ in the resolvent; we have 
only left out the principal symbol of $G_{\lambda }$. (It is not
clear what kind of operator comes out of applying the log Cauchy
formula to this term in general.)

\proclaim{Theorem 5.2} Consider a normal elliptic
realization $(P_++G)_T$, where
$P$ is integer order $m>0$, $G$ is of order and class $m$, 
and $T=\{T_0,\dots,T_{m-1}\}$ is  normal,
with 
entries $T_k$ of order and class $k$. Assume that $m>n$. 

With $B$ defined in Lemma {\rm 5.1}, we have
that $$
C_0(I, (P_++G)_T )= -\tfrac 1m\operatorname{res} ((\log P
)_+)-\tfrac1m\operatorname{res} (B).\tag5.9
$$
Here $$c^P_{n,+}=-\tfrac 1m(\operatorname{res} ((\log P
)_+), \quad
c^G_n=-\tfrac 1m\operatorname{res}(B).\tag5.10
$$ 
\endproclaim 

\demo{Proof} This goes as in Theorems  3.9 and 4.9. The necessary
symbol information has been provided above, so we just have to
identify the contributions from the specific homogeneous terms.\qed
\enddemo

In some cases one can get a more informative formula, as 
the following example (similar to \cite{GSc2, Rem.\ 4.2}) shows.
 
 \example{Example 5.3}
Consider a second-order strongly elliptic differential operator $P$,
of the form $$
P=-\partial _{x_n}^2+P'\tag5.11
$$
in a collar neighborhood of $X'$, where $P'$ is a positive
selfadjoint second-order elliptic operator on $X'$. Let $T=\gamma
_0$, restriction to $X'$; then $(P_+)_{\gamma _0}$ is the Dirichlet
realization of $P$. The resolvent $R_\lambda $ does not have high
enough order to be trace-class, but we can iterate it, considering
$$
R_\lambda ^N=\tfrac{\partial _\lambda ^{N-1}}{(N-1)!}R_\lambda
=\tfrac{\partial _\lambda ^{N-1}}{(N-1)!}Q_{\lambda ,+}+\tfrac{\partial
_\lambda ^{N-1}}{(N-1)!}G_\lambda   =(Q_\lambda )^N_++G_\lambda ^{(N)}\tag5.12
$$ for
$N>n/2$ instead. It is easily verified (further details in \cite{GSc2,
Rem.\ 4.2}) that $\tr_nG_\lambda =-\frac14(P'-\lambda )^{-1}$, a
resolvent on $X'$ (times a constant). The interior contribution to
the coefficient of $(-\lambda )^{-N}$ is $$
-\tfrac12 \operatorname{res}_X((\log P)_+),\tag5.13
$$
in view of the considerations in Remark 3.12. The same considerations
plus the information from Section 1 for closed manifolds, applied to
$-\frac14(P'-\lambda )^{-1}$, gives that the s.g.o.\ 
contribution is$$
\tfrac18 \operatorname{res}_{X'}(\log P').\tag5.14
$$ 
So here$$
C_0(I,(P_+)_{\gamma _0})=-\tfrac12 \operatorname{res}_X((\log P)_+)+\tfrac18 \operatorname{res}_{X'}(\log P'),\tag5.15
$$
where we have logarithmic operators in both terms.

It may be remarked as in \cite{GSc2} that the interior term vanishes
when $n$ is odd, the boundary term vanishes when $n$ is even.

\endexample

\comment
[Den anden trace defect: Man kan komme s\aa{} langt som at skrive en
f-agtig formel op, og g\o{}re det til et integral med log, men
interpretationen? Det er bare $G^\pm (\log P_1)$'s principale led,
der ikke er d\ae{}kning for, det m\aa{} afvente yderligere teori (og
er m\aa{}ske begr\ae{}nset interessant. --- eller dukker der en
rigtig kommutator op?]

\subhead 4. The commutator trace defect \endsubhead

When one tries to generalize the trace defect formula (2.6) to
$\psi $dbo's, one meets the problem that
singular Green operators such as $G^\pm(\log P)$ are not in the
calculus, and they appear in non-commuting situations here, e.g.\ 
in calculations such as $$
\multline
P_+(\log P_1)_+-
(\log P_1)_+P_+\\
=(P\log P_1-\log P_1 P)_+ -G^+(P)G^-(\log
P_1)+G^+(\log P_1)G^-(P) .\endmultline\tag4.1
$$ 
Also the lack of commuativity of $\psi $do's with s.g.o.s gives trouble,
e.g.\ in the study of $G(\log P_1)_+-(\log P_1)_+G$ and its
compositions with other operators.

One possibility might be that one could give a residue formula involving
only the classical part of $\log P_1$. Consider again its definition.
When the principal symbol is scalar (as we assume here),
the terms in the symbol (cf.\ (1.14)ff.) are calculated by:
$$\aligned
\tfrac i{2\pi }\int_{\Cal C'}\log\lambda
\,(p_{1,m}(x,\xi )-\lambda )^{-1}I\, d\lambda &=\log (p_{1,m}(x,\xi
))\,I\\
&=m\log (|\xi 
|)\,I+\log (p_{1,m}(x,\xi /|\xi |))\,I\text{ for }|\xi |\ge 1,\\
\tfrac i{2\pi }\int_{\Cal C'}\log\lambda
\,q_{-m-j}(x,\xi ,\lambda )\, d\lambda &=b_{-j}(x,\xi ), \text{ all }j\ge 1,
\endaligned\tag 4.2$$
where in the last line, the fact that $q_{-m-j}$ consists of terms
containing 
$(p_{1,m}-\lambda )^{-1}$ in powers $\ge 2$ allows integrations by
part replacing $\log \lambda $ by $\lambda ^{-1}$, leading to a
standard homogeneous symbol of degree $-j$. So
$$
\log P=\op(m\log [\xi ]I)+B,\quad B=\op(b),
\tag4.3$$ 
with symbol $b\sim \sum_{j\ge 0}b_{-j}$ having$$
b_0(x,\xi )=\log(p_{1,m}(x,\xi /|\xi|))\,I \text{ for }|\xi |\ge 1\tag 4.4
$$
(clearly homogeneous of degree 0); the other $b_{-j}$ being defined in (4.2).

It may be so, that the desired formulas will be valid with $B$ instead of
$\log P_1$ itself, or even with $B$ replaced by
$$
B'=\op(b'),\quad b'=b-b_0.\tag4.4
$$

$B$ or $B'$ can be viewed as ``the log without the log'', moving
us still further away from a genuine use of an operator-logarithm here.

Let ... As in (2.21), we have that $\Tr T=\Tr R$, where 
$$
T=[A,A']Q_{\lambda ,+},\quad R=A[A',Q_{\lambda,+ }] ,
\tag4.7
$$
and we shall work on the formulas for $R$. We have to deal with three
types of cases:
$$\aligned
\text{Case 1: }A&=P_+,\quad A'=P_+,\\ 
R&=P_+[P'_+,Q_{\lambda ,+}]=
\endaligned$$

[Og s\aa{} ser regningerne alligevel uoverkommelige ud, og det er
m\aa{}ske en falsk conjecture at kun $B$ skulle indg\aa{} --- dog er
det faktisk tilf\ae{}ldet for den f\o{}rste trace defekt formel.]
\endcomment

\enddocument